\documentclass[a4paper,11pt]{article}
\usepackage{amsfonts,amssymb,amsmath,epsf,eepic,graphics,graphicx,bbm}

\usepackage[english]{babel}
 \makeatletter
 \@addtoreset{enunciato}{section}
 \makeatother

 \newcounter{enunciato}[section]

 \newtheorem{ittheorem}{Theorem}
 \newtheorem{itlemma}{Lemma}
 \newtheorem{itproposition}{Proposition}
 \newtheorem{itdefinition}{Definition}
 \newtheorem{itclaim}{Claim}
 \newtheorem{itfact}{Fact}
 \newtheorem{itconjecture}{Conjecture}

 \newenvironment{theorem}{\addtocounter{enunciato}{1}
 \begin{ittheorem}}{\end{ittheorem}}

 \newenvironment{lemma}{\addtocounter{enunciato}{1}
 \begin{itlemma}}{\end{itlemma}}

 \newenvironment{proposition}{\addtocounter{enunciato}{1}
 \begin{itproposition}}{\end{itproposition}}

 \newenvironment{definition}{\addtocounter{enunciato}{1}
 \begin{itdefinition}}{\end{itdefinition}}

 \newenvironment{claim}{\addtocounter{enunciato}{1}
 \begin{itclaim}}{\end{itclaim}}

 \newenvironment{fact}{\addtocounter{enunciato}{1}
 \begin{itfact}}{\end{itfact}}

 \newenvironment{conjecture}{\addtocounter{enunciato}{1}
 \begin{itconjecture}}{\end{itconjecture}}

 \newcommand{\be}[1]{\begin{equation}\label{#1}}
 \newcommand{\ee}{\end{equation}}

 \newcommand{\bea}[1]{\begin{eqnarray}\label{#1}}
 \newcommand{\eea}{\end{eqnarray}}

 \newcommand{\bl}[1]{\begin{lemma}\label{#1}}
 \newcommand{\el}{\end{lemma}}

 \newcommand{\bt}[1]{\begin{theorem}\label{#1}}
 \newcommand{\et}{\end{theorem}}

 \newcommand{\bd}[1]{\begin{definition}\label{#1}}
 \newcommand{\ed}{\end{definition}}

 \newcommand{\bcl}[1]{\begin{claim}\label{#1}}
 \newcommand{\ecl}{\end{claim}}

 \newcommand{\bfact}[1]{\begin{fact}\label{#1}}
 \newcommand{\efact}{\end{fact}}

 \newcommand{\bp}[1]{\begin{proposition}\label{#1}}
 \newcommand{\ep}{\end{proposition}}

 \newcommand{\bc}[1]{\begin{corollary}\label{#1}}
 \newcommand{\ec}{\end{corollary}}

 \newcommand{\bcj}[1]{\begin{conjecture}\label{#1}}
 \newcommand{\ecj}{\end{conjecture}}

 \newcommand{\bpr}{\begin{proof}}
 \newcommand{\epr}{\end{proof}}

 \newcommand{\bprl}[1]{\begin{proofof}{\it\ref{#1}}.\,\,}
 \newcommand{\eprl}{\end{proofof}}

 \newcommand{\bi}{\begin{itemize}}
 \newcommand{\ei}{\end{itemize}}

 \newcommand{\ben}{\begin{enumerate}}
 \newcommand{\een}{\end{enumerate}}

\newenvironment{proof}{\noindent {\em Proof}.\,\,}
{\hspace*{\fill}$\halmos$\medskip}
\newenvironment{proofof}{\noindent {\em Proof of Lemma\,\,}}
{\hspace*{\fill}$\halmos$\medskip}
\newcommand{\halmos}{\rule{1ex}{1.4ex}}

\def \Z {{\mathbb Z}}

\def \N {{\mathbb N}}

\def \ba {\begin{array}}
\def \ea {\end{array}}

\def \P {{\mathbb P}}
\def \E {{\mathbb E}}
\def \U {{\mathbb U}}

\def \cX {{\mathcal X}}

\def \cC {{\mathcal C}}
\def \cP {{\mathcal P}}
\def \cS {{\mathcal S}}
\def \cD {{\mathcal D}}
\def \cE {{\mathcal E}}
\def \cA {{\mathcal A}}
\def \cB {{\mathcal B}}

\def \cW {{\mathcal W}}
\def \cI {{\mathcal W}}
\def \cU {{\mathcal U}}
\def \cI {{\mathcal I}}
\def \cV {{\mathcal V}}

\def \L {\Lambda}
\def \CAPA {{\hbox{\footnotesize\rm CAP}}}

\def\b{\beta}
\def\d{\delta}
\def\g{\gamma}

\def\s{\sigma}
\def\t{\tau}

\def\D{\Delta}
\def\L{\Lambda}
\def\G{\Gamma}

\def\e{\epsilon}
\def\1{\mathbbm{1}}

\begin{document}


\title{Homogeneous nucleation for Glauber and Kawasaki dynamics\\
in large volumes at low temperatures}

\author{
\renewcommand{\thefootnote}{\arabic{footnote}}
A.\ Bovier
\footnotemark[1]\,\,\,\footnotemark[2]
\\
\renewcommand{\thefootnote}{\arabic{footnote}}
F.\ den Hollander
\footnotemark[3]\,\,\,\footnotemark[4]
\\
\renewcommand{\thefootnote}{\arabic{footnote}}
C.\ Spitoni
\footnotemark[3]\,\,\,\footnotemark[4]
}

\footnotetext[1]{
Weierstrass-Institut f\"ur Angewandte Analysis und Stochastik,
Mohrenstrasse 39, 10117 Berlin, Germany
}
\footnotetext[2]{
Institut f\"ur Mathematik, Technische Universit\"at Berlin,
Strasse des 17.\ Juni 136, 10623 Berlin, Germany
}
\footnotetext[3]{
Mathematical Institute, Leiden University, P.O.\ Box 9512,
2300 RA Leiden, The Netherlands
}
\footnotetext[4]{
EURANDOM, P.O.\ Box 513, 5600 MB Eindhoven, The Netherlands
}

\maketitle

\begin{abstract}

In this paper we study metastability in large volumes at low temperatures.
We consider both Ising spins subject to Glauber spin-flip dynamics and lattice
gas particles subject to Kawasaki hopping dynamics. Let $\b$ denote the
inverse temperature and let $\L_\b \subset \Z^2$ be a square box with
periodic boundary conditions such that $\lim_{\b\to\infty}|\L_\b|=\infty$.
We run the dynamics on $\L_\b$ starting from a random initial configuration where
all the droplets (= clusters of plus-spins, respectively, clusters of particles)
are small. For large $\b$, and for interaction parameters that correspond to the
metastable regime, we investigate how the transition from the metastable state (with
only small droplets) to the stable state (with one or more large droplets) takes
place under the dynamics. This transition is triggered by the appearance of a
single \emph{critical droplet} somewhere in $\L_\b$. Using potential-theoretic
methods, we compute the \emph{average nucleation time} (= the first time a
critical droplet appears and starts growing) up to a multiplicative factor
that tends to one as $\b\to\infty$. It turns out that this time grows as
$Ke^{\Gamma\b}/|\L_\b|$ for Glauber dynamics and $K\b e^{\Gamma\b}/|\L_\b|$
for Kawasaki dynamics, where $\Gamma$ is the local canonical, respectively,
grand-canonical energy to create a critical droplet and $K$ is a constant
reflecting the geometry of the critical droplet, provided these times tend to
infinity (which puts a growth restriction on $|\L_\b|$). The fact that the average
nucleation time is inversely proportional to $|\L_\b|$ is referred to as
\emph{homogeneous nucleation}, because it says that the critical droplet for
the transition appears essentially independently in small boxes that partition
$\L_\b$.

{\it Acknowledgment.} The authors thank Alessandra Bianchi, Alex Gaudilli\`ere,
Dima Ioffe, Francesca Nardi, Enzo Olivieri and Elisabetta Scoppola for ongoing
discussions on meta\-stability and for sharing their work in progress. CS thanks
Martin Slowik for stimulating exchange. AB and FdH are supported by DFG and NWO
through the Dutch-German Bilateral Research Group on ``Random Spatial Models from
Physics and Biology'' (2003--2009). CS is supported by NWO through grant 613.000.556.

\end{abstract}

\newpage


\section{Introduction and main results}
\label{S1}

\subsection{Background}
\label{Intro}

In a recent series of papers, Gaudilli\`ere, den Hollander, Nardi, Olivieri, and
Scoppola~\cite{GdHNOS07a,GdHNOS07b,GdHNOS07c} study a system of lattice gas particles
subject to Kawasaki hopping dynamics in a \emph{large} box at low temperature and low
density. Using the so-called \emph{path-wise approach} to metastability (see Olivieri
and Vares~\cite{OV04}), they show that the transition time between the metastable state
(= the gas phase with only small droplets) and the stable state (= the liquid phase with
one or more large droplets) is inversely proportional to the volume of the large box,
provided the latter does not grow too fast with the inverse temperature. This type of
behavior is called \emph{homogeneous nucleation}, because it corresponds to the situation
where the critical droplet triggering the nucleation appears essentially independently
in small boxes that partition the large box. The \emph{nucleation time} (= the first
time a critical droplet appears and starts growing) is computed up to a multiplicative
error that is small on the scale of the \emph{exponential} of the inverse temperature.
The techniques developed in \cite{GdHNOS07a,GdHNOS07b,GdHNOS07c} center around the
idea of \emph{approximating} the low temperature and low density Kawasaki lattice gas
by an \emph{ideal gas} without interaction and showing that this ideal gas stays \emph{close
to equilibrium} while exchanging particles with droplets that are growing and shrinking.
In this way, the large system is shown to behave essentially like the union of many small
independent systems, leading to homogeneous nucleation. The proofs are long and complicated,
but they provide considerable detail about the \emph{typical trajectory} of the system
prior to and shortly after the onset of nucleation.

In the present paper we consider the same problem, both for Ising spins subject to Glauber
spin-flip dynamics and for lattice gas particles subject to Kawasaki hopping dynamics.
Using the \emph{potential-theoretic approach} to metastability (see Bovier~\cite{B07}),
we improve \emph{part} of the results in \cite{GdHNOS07a,GdHNOS07b,GdHNOS07c}, namely,
we compute the \emph{average} nucleation time up to a multiplicative error that \emph{tends
to one} as the temperature tends to zero, thereby providing a very sharp estimate of the
time at which the gas starts to condensate.

We have no results about the typical time it takes for the system to grow a large droplet
after the onset of nucleation. This is a hard problem that will be addressed in future work.
All that we can prove is that the dynamics has a negligible probability to shrink down a
supercritical droplet once it has managed to create one. At least this shows that the
appearance of a single critical droplet indeed represents the threshold for nucleation, as
was shown in \cite{GdHNOS07a,GdHNOS07b,GdHNOS07c}. A further restriction is that we need to
draw the initial configuration according to a \emph{class} of initial distributions on the
set of subcritical configurations, called the last-exit biased distributions, since these are
particularly suitable for the use of potential theory. It remains a challenge to investigate
to what extent this restriction can be relaxed. This problem is addressed with some success
in \cite{GdHNOS07a,GdHNOS07b,GdHNOS07c}, and will also be tackled in future work.

Our results are an extension to large volumes of the results for small volumes
obtained in Bovier and Manzo~\cite{BM02}, respectively, Bovier, den Hollander, and
Nardi~\cite{BdHN06}. In large volumes, even at low temperatures \emph{entropy}
is competing with energy, because the metastable state and the states that evolve
from it under the dynamics have a highly non-trivial structure. Our main goal in
the present paper is to extend the potential-theoretic approach to metastability
in order to be able to deal with large volumes. This is part of a broader programme
where the objective is to adapt the potential-theoretic approach to situations where
entropy cannot be neglected. In the same direction, Bianchi, Bovier, and Ioffe~\cite{BBI}
study the dynamics of the random field Curie-Weiss model on a finite box at a fixed
\emph{positive} temperature.

As we will see, the basic difficulty in estimating the nucleation time is to obtain sharp
upper and lower bounds on \emph{capacities}. Upper bounds follow from the Dirichlet
variational principle, which represents a capacity as an infimum over a class of test
functions. In \cite{BBI} a new technique is developed, based on a variational principle
due to Berman and Konsowa~\cite{BK90}, which represent a capacity as a supremum over a
class of unit flows. This technique allows for getting lower bounds and it will be exploited
here too.

\subsection{Ising spins subject to Glauber dynamics}
\label{S1.1}

We will study models in finite boxes, $\L_\b$, in the limit as both the inverse temperature,
$\b$, and the volume of the box, $|\L_\b|$, tend to infinity. Specifically, we let
$\L_\b \subset \Z^2$ be a square box with odd side length, centered at the origin
with periodic boundary conditions.  A spin configuration is denoted by $\s=\{\s(x)
\colon\,x\in\L_\b\}$, with $\s(x)$ representing the spin at site $x$, and is an element
of $\cX_\b = \{-1,+1\}^{\L_\b}$. It will frequently be convenient to identify a configuration
$\s$ with its \emph{support},  defined as $\textnormal{supp}[\s]=\{x\in\L_\b\colon\,\s(x)=+1\}$.

The interaction is defined by the the usual Ising Hamiltonian
\be{Ham}
H_\b(\s) = - \frac{J}{2} \sum_{ {(x,y) \in \L_\b} \atop {x \sim y} }
\s(x)\s(y) - \frac{h}{2} \sum_{x \in \L_\b} \s(x),
\qquad \s\in\cX_\b,
\ee
where $J>0$ is the pair potential, $h>0$ is the magnetic field, and $x \sim y$ means that
$x$ and $y$ are nearest neighbors. The \emph{Gibbs measure} associated with $H_\b$ is
\be{Gibbs}
\mu_\b(\s) = \frac{1}{Z_\b}\,e^{-\b H_\b(\s)}, \qquad \s \in \cX_\b,
\ee
where $Z_\b$ is the normalizing partition function.

The dynamics of the model will the a continuous-time Markov chain, $(\s(t))_{t\geq 0}$,
with state space $\cX_\b$ whose transition rates are given by
\be{rate}
c_\b(\s,\s') = \left\{
\begin{array}{ll}
e^{-\b [H_\b(\s')-H_\b(\s)]_+},
&\mbox{for } \s'=\s^x \mbox{ for some } \,x \in \L_\b,\\
0,
&\mbox{otherwise},
\end{array}
\right.
 \ee
where $\s^x$ is the configuration obtained from $\s$ by flipping the spin at site $x$. We
refer to this Markov process as \emph{Glauber dynamics}. It is ergodic and reversible with
respect to its unique invariant measure, $\mu_\b$, i.e.,
\be{rev}
\mu_\b(\s) c_\b(\s,\s') = \mu_\b(\s') c_\b(\s',\s),
\qquad \forall\,\s,\s' \in \cX_\b.
 \ee

\begin{figure}[htbp]
\centering
\includegraphics[height=4cm]{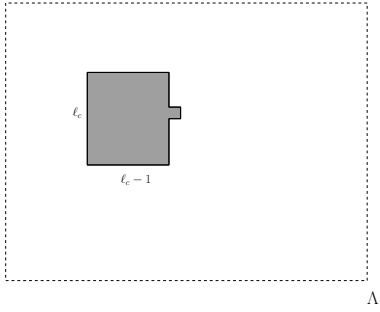}
\caption{\small
A critical droplet for Glauber dynamics on $\L$. The shaded
area represents the $(+1)$-spins, the non-shaded area the $(-1)$-spins
(see (\ref{lcdef})).}
\label{fig:critdropG}
\end{figure}

\begin{figure}[htbp]
\vspace{0.5cm}
\centering
\includegraphics[height=4cm]{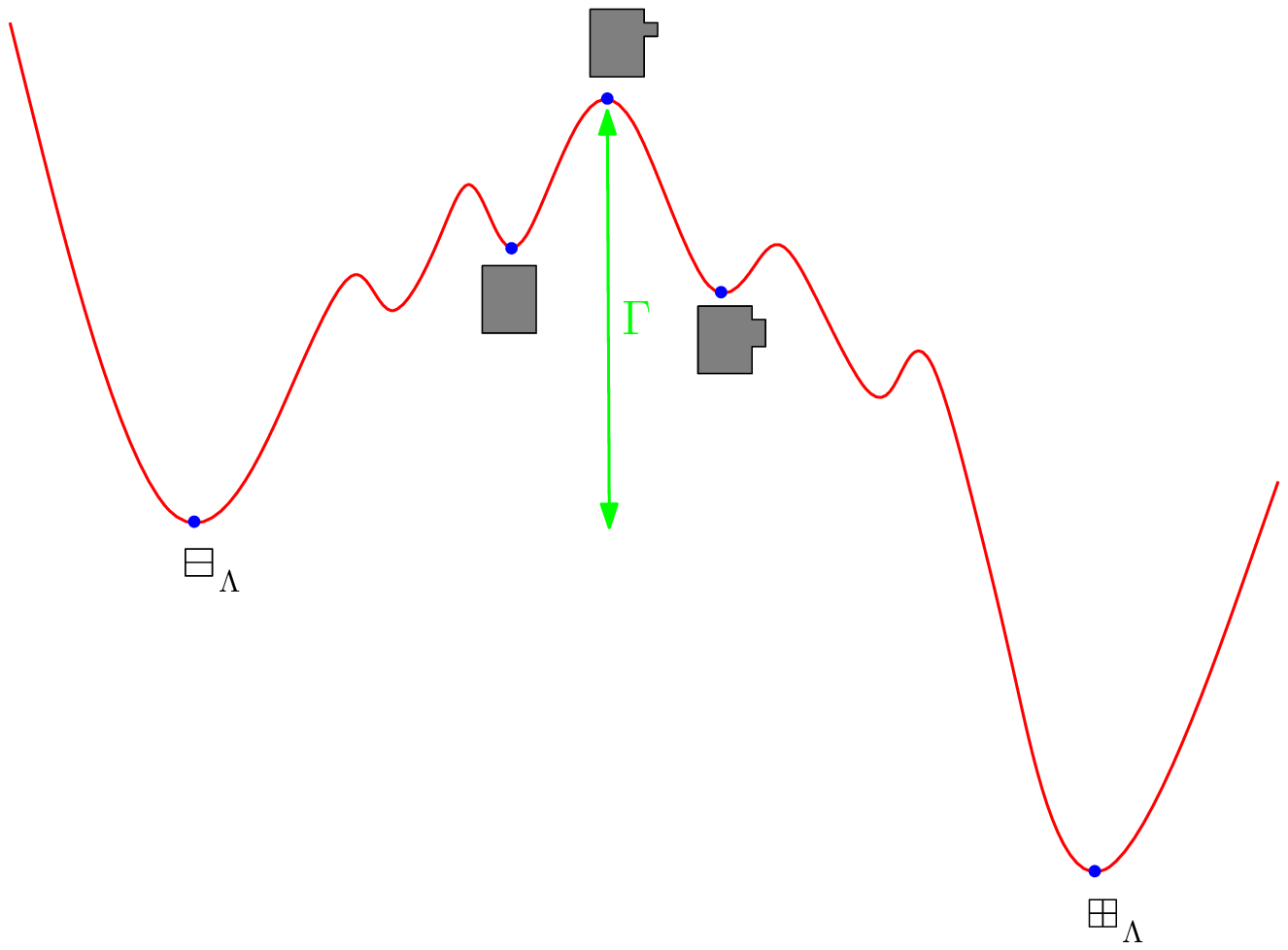}
\caption{\small
A nucleation path from $\boxminus_{\L}$ to $\boxplus_{\L}$ for Glauber dynamics.
$\Gamma$ in (\ref{Gamdef}) is the minimal energy barrier the path has to overcome
under the local variant of the Hamiltonian in (\ref{Ham}).}
\label{fig:nuclpath}
\end{figure}

Glauber dynamics exhibits metastable behavior in the regime
\be{metreg}
0 < h < 2J, \qquad \b \to \infty.
\ee
To understand this, let us briefly recall what happens in a finite $\b$-independent
box $\L\subset\Z^2$. Let $\boxminus_{\L}$ and $\boxplus_{\L}$ denote the configurations
where all spins in $\L$ are $-1$, respectively, $+1$. As was shown by Neves and
Schonmann~\cite{NS91}, for Glauber dynamics \emph{restricted} to $\L$ with periodic
boundary conditions and subject to (\ref{metreg}), the \emph{critical droplets} for
the crossover from $\boxminus_{\L}$ to $\boxplus_{\L}$ are the set of all those configurations
where the $(+1)$-spins form an $\ell_c\times(\ell_c-1)$ quasi-square (in either of both
orientations) with a protuberance attached to one of its longest sides, where
\be{lcdef}
\ell_c = \left\lceil\frac{2J}{h}\right\rceil
\ee
(see Figs.~\ref{fig:critdropG} and \ref{fig:nuclpath}; for non-degeneracy reasons
it is assumed that $2J/h\notin\N$). The quasi-squares without the protuberance are
called \emph{proto-critical droplets}.

Let us now return to our setting with finite $\b$-dependent volumes $\L_\b
\subset\Z^2$. We will start our dynamics on $\L_\b$ from initial configurations
in which all droplets are ``sufficiently small''. To make this notion
precise, let $C_B(\s)$, $\s\in\cX_\b$, be the configuration that
is obtained from $\s$ by a ``bootstrap percolation map'', i.e., by
circumscribing all the droplets in $\s$ with rectangles, and continuing
to doing so in an iterative manner until a union of disjoint rectangles
is obtained (see Koteck\'y and Olivieri~\cite{KO93}). We call $C_B(\s)$
\emph{subcritical} if all its rectangles fit inside a proto-critical
droplet and are at distance $\geq 2$ from each other (i.e., are
non-interacting).

\bd{Rdef}
(a) $\cS=\{\s\in\cX_\b\colon\,C_B(\s) \mbox{ is subcritical}\,\}$.\\
(b) $\cP=\{\s\in\cS\colon c_\b(\s,\s')>0 \mbox{ for some }
\s'\in
\cS^c\}$.\\
(c) $\cC=\{\s'\in\cS^c\colon c_\b(\s,\s')>0 \mbox{ for
  some } \s\in \cS\}$.
\ed

\noindent
We refer to $\cS$, $\cP$ and $\cC$ as the set of  \emph{subcritical},
\emph{proto-critical}, respectively, \emph{critical} configurations. Note
that, for ever $\s\in\cX_\b$, each step in the bootstrap percolation map
$\s\to C_B(\s)$ deceases the energy, and therefore the Glauber dynamics
moves from $\s$ to $C_B(\s)$ in a time of order one. This is why $C_B(\s)$
rather than $\s$ appears in the definition of $\cS$.

For $\ell_1,\ell_2\in\N$, let $R_{\ell_1,\ell_2}(x) \subset \L_\b$ be
the $\ell_1\times\ell_2$ rectangle whose lower-left corner is $x$. We
always take $\ell_1\leq\ell_2$ and allow for both orientations of the
rectangle. For $L=1,\dots,2\ell_c-3$, let $Q_L(x)$ denote the $L$-th
element in the \emph{canonical} sequence of growing squares and quasi-squares
\begin{equation}
\label{qsqdef}
R_{1,2}(x),\,R_{2,2}(x),\,R_{2,3}(x),\,R_{3,3}(x),\dots,\,
R_{\ell_c-1,\ell_c-1}(x),\,R_{\ell_c-1,\ell_c}(x).
\end{equation}

In what follows we will choose to start the dynamics in a way that is
suitable for the use of potential theory, as follows. First, we take
the initial law to be concentrated on sets $S_L\subset\cS$ defined by
\be{SLdef}
S_L = \left\{\s\in\cS\colon\,\mbox{ each rectangle in } C_B(\s) \mbox{ fits
inside } Q_L(x) \mbox{ for some } x\in\L_\b\right\},
\ee
where $L$ is any integer satisfying
\be{L*def}
L^* \leq L \leq 2\ell_c-3 \quad \mbox{ with } \quad
L^* = \min \left\{1 \leq L \leq 2\ell_c-3\colon\,
\lim_{\b\to\infty} \frac{\mu_\b(\cS_L)}{\mu_\b(\cS)} = 1\right\}.
\ee
In words, $\cS_L$ is the subset of those subcritical configurations whose droplets
fit inside a square or quasi-square labeled $L$, with $L$ chosen large enough so
that $\cS_L$ is typical within $\cS$ under the Gibbs measure $\mu_\b$ as $\b\to\infty$
(our results will not depend on the choice of $L$ subject to these restrictions). Second,
we take the initial law to be biased according to the last exit of $\cS_L$ for the
transition from $\cS_L$ to a target set in $\cS^c$. (Different choices will be made
for the target set, and the precise definition of the biased law will be given in
Section~\ref{S2.1}.) This is a \emph{highly specific} choice, but clearly one of
physical interest.

\medskip\noindent
{\bf Remarks:}
(1) Note that $\cS_{2\ell_c-3}=\cS$, which implies that the range of $L$-values in
(\ref{L*def}) is non-empty. The value of $L^*$ depends on how fast $\L_\b$ grows with
$\b$. In Appendix~\ref{appC.1} we will show that, for every $1 \leq L \leq 2\ell_c-4$,
$\lim_{\b\to\infty}\mu_\b(\cS_L)/\mu_\b(\cS)=1$ if and only if $\lim_{\b\to\infty}
|\L_\b|e^{-\b\Gamma_{L+1}}=0$ with $\Gamma_{L+1}$ the energy needed to create a droplet
$Q_{L+1}(0)$ at the origin. Thus, if $|\L_\b|=e^{\theta\beta}$, then $L^*= L^*(\theta)
= (2\ell_c-3) \wedge \min\{L\in\N\colon\,\Gamma_{L+1}>\theta\}$, which increases stepwise
from $1$ to $2\ell_c-3$ as $\theta$ increases from $0$ to $\Gamma$ defined in
(\ref{Gamdef}).\\
(2) If we draw the initial configuration $\s_0$ from some subset of $\cS$ that
has a strong recurrence property under the dynamics, then the choice of initial
distribution on this subset should not matter. This issue will be addressed in
future work.

\medskip
To state our main theorem for Glauber dynamics, we need some further notation.
The key quantity for the nucleation process is
\be{Gamdef}
\Gamma = J[4\ell_c]-h[\ell_c(\ell_c-1)+1],
\ee
which is the energy needed to create a critical droplet of $(+1)$-spins at a given
location in a sea of $(-1)$-spins (see Figs.~\ref{fig:critdropG} and \ref{fig:nuclpath}).
For $\s\in\cX_\b$, let $\P_\s$ denote the law of the dynamics starting from $\s$ and,
for $\nu$ a probability distribution on $\cX$, put
\be{PGldef}
\P_\nu(\cdot) = \sum_{\s\in\cX_\b} \P_\s(\cdot)\,\nu(\s).
\ee
For a non-empty set $\cA\subset\cX_\b$, let
\be{tauAdef}
\tau_\cA = \inf\{t>0\colon\,\s_t\in\cA,\,\s_{t^-}\notin\cA\}
\ee
denote the first time the dynamics enters $\cA$. For non-empty and disjoint sets
$\cA,\cB\subset\cX_\b$, let $\nu_\cA^\cB$ denote the \emph{last-exit biased
distribution on $\cA$ for the crossover to $\cB$} defined in (\ref{lastexitdef})
in Section~\ref{S2.1}. Put
\be{N12def}
N_1=4\ell_c, \qquad N_2=\tfrac43(2\ell_c-1).
\ee
For $M\in\N$ with $M \geq \ell_c$, define
\be{DMdef}
\cD_M = \big\{\s\in\cX_\b\colon\,\exists\,x\in\L_\b \mbox{ such that }
\textnormal{supp}[C_B(\s)] \supset R_{M,M}(x)\big\},
\ee
i.e., the set of configurations containing a supercritical droplet of size $M$. For
our results below to be valid we need to assume that
\be{Lambdaconds}
\lim_{\b\to\infty} |\L_\b|=\infty,\qquad \lim_{\beta\to\infty} |\L_\b|\,e^{-\b\G}=0.
\ee

\bt{HNGlauber}
In the regime {\rm (\ref{metreg})}, subject to {\rm (\ref{L*def}) and
(\ref{Lambdaconds})}, the following hold:\\
(a)
\be{crossbdsGlauber}
\lim_{\b\to\infty} |\L_\b|\,
e^{-\b\Gamma}\,\E_{\nu_{\cS_L}^{\cS^c}}\left(\tau_{\cS^c}\right) = \frac{1}{N_1}.
\ee
(b)
\be{crossbdsGlauber2}
\lim_{\b\to\infty} |\L_\b|\,
e^{-\b\Gamma}\,\E_{\nu_{\cS_L}^{\cS^c\backslash\cC}}\left(\tau_{\cS^c\backslash\cC}\right)
= \frac{1}{N_2}.
\ee
(c)
\be{crossbdsGlauber3}
\lim_{\b\to\infty} |\L_\b|\, e^{-\b\Gamma}\,
\E_{\nu_{\cS_L}^{\cD_M}}\left(\tau_{\cD_M}\right) = \frac{1}{N_2},
\qquad \forall\,\ell_c \leq M \leq 2\ell_c-1.
\ee
\et

\medskip\noindent
The proof of Theorem~\ref{HNGlauber} will be given in Section~\ref{S3}.
Part (a) says that the average time to create a critical droplet is $[1+o(1)]
e^{\b\Gamma}/N_1|\L_\b|$. Parts (b) and (c) say that the average time to go
beyond this critical droplet and to grow a droplet that is twice as large is
$[1+o(1)]e^{\b\Gamma}/N_2|\L_\b|$. The factor $N_1$ counts the number of
\emph{shapes} of the critical droplet, while $|\L_\b|$ counts the number of
\emph{locations}. The average times to create a critical, respectively, a
supercritical droplet differ by a factor $N_2/N_1<1$. This is because once
the dynamics is ``on top of the hill'' $\cC$ it has a positive probability to
``fall back'' to $\cS$. On average the dynamics makes $N_1/N_2>1$ attempts to
reach the top $\cC$ before it finally ``falls over'' to $\cS^c\backslash\cC$.
After that, it rapidly grows a large droplet (see Fig.~\ref{fig:nuclpath}).

\medskip\noindent
{\bf Remarks:}
(1) The second condition in (\ref{Lambdaconds}) will not actually be used in the
proof of Theorem~\ref{HNGlauber}(a). If this condition fails, then there is a positive
probability to see a proto-critical droplet in $\L_\b$ under the starting measure
$\nu_{\cS_L}^{\cS^c}$, and nucleation sets in immediately. Theorem~\ref{HNGlauber}(a)
continues to be true, but it no longer describes metastable behavior.\\
(2) In Appendix~\ref{appD} we will show that the \emph{average} probability under
the Gibbs measure $\mu_\b$ of destroying a supercritical droplet and returning to
a configuration in $\cS_L$ is exponentially small in $\b$. Hence, the crossover from
$\cS_L$ to $\cS^c\backslash\cC$ represents the true threshold for nucleation, and
Theorem~\ref{HNGlauber}(b) represents the true nucleation time.\\
(3) We expect  Theorem~\ref{HNGlauber}(c) to hold for values of $M$ that grow with $\b$
as $M=e^{o(\b)}$. As we will see in Section~\ref{S3.3}, the necessary capacity
estimates carry over, but the necessary equilibrium potential estimates are not
yet available. This problem will be addressed in future work.\\
(4) Theorem~\ref{HNGlauber} should be compared with the results in Bovier and
Manzo~\cite{BM02} for the case of a finite $\b$-independent box $\L$ (large
enough to accommodate a critical droplet). In that case, if the dynamics starts
from $\boxminus_{\L}$, then the average time it needs to hit $\cC_{\L}$ (= the
set of configurations in $\L$ with a critical droplet), respectively, $\boxplus_{\L}$
equals
\be{Kasymp}
Ke^{\b\Gamma}[1+o(1)], \mbox{ with } K = K(\L,\ell_c)
= \frac{1}{N}\,\frac{1}{|\L|} \mbox{ for } N=N_1,N_2.
\ee
(4) Note that in Theorem~\ref{HNGlauber} we compute the first time when a critical
droplet appears anywhere (!) in the box $\L_\b$. It is a different issue to compute
the first time when the plus-phase appears near the origin. This time, which depends
on how a supercritical droplet grows and eventually invades the origin, was studied
by Dehghanpour and Schonmann~\cite{DS97a,DS97b}, Shlosman and Schonmann~\cite{SS98}
and, more recently, by Cerf and Manzo~\cite{CM08}.

\subsection{Lattice gas subject to Kawasaki dynamics}
\label{S1.2}

We next consider the lattice gas subject to Kawasaki dynamics and state a similar
result for homogeneous nucleation. Some aspects are similar as for Glauber dynamics,
but there are notable differences.

A lattice gas configuration is denoted by $\s=\{\s(x)\colon\,x\in\cX_\b\}$, with $\s(x)$
representing the number of particles at site $x$, and is an element of
$\cX_\b = \{0,1\}^{\L_\b}$. The Hamiltonian is given by
\be{Hamc}
H_\b(\s) = - U \sum_{ {(x,y) \in \L_\b} \atop {x \sim y} }
\s(x)\s(y), \qquad \s\in\cX_\b,
\ee
where $-U<0$ is the binding energy and $x \sim y$ means that $x$ and $y$ are neighboring
sites. Thus, we are working in the \emph{canonical ensemble}, i.e., there is no term analogous
to the second term in (\ref{Ham}). The  number of particles in $\L_\b$ is
\be{nopart}
n_\b = \lceil\,\rho_\b|\L_\b|\,\rceil,
\ee
where $\rho_\b$ is the particle density, which is chosen to be
\be{rhodef}
\rho_\b = e^{-\b\Delta}, \qquad \Delta>0.
\ee
Put
\be{cXbetandef}
\cX_\b^{(n_\b)} = \{\s\in\cX_\b\colon\,|\textnormal{supp}[\s]|=n_\b\},
\ee
where $\textnormal{supp}[\s]=\{x\in\L_\b\colon\,\s(x)=1\}$.

\medskip\noindent
{\bf Remark:}
If we were to work in the \emph{grand-canonical ensemble}, then we would have to
consider the Hamiltonian
\be{Hamgc}
H^{gc}(\s) = - U \sum_{ {(x,y) \in \L_\b} \atop {x \sim y} } \s(x)\s(y)
+ \Delta \sum_{x\in\L_\b} \s(x), \qquad \s\in\cX_\b,
\ee
with $\Delta>0$ an activity parameter taking over the role of $h$ in (\ref{Ham}).
The second term would mimic the presence of an infinite gas reservoir with
density $\rho_\b$ outside $\L_\b$. Such a Hamiltonian was used in earlier work on
Kawasaki dynamics, when a finite $\b$-independent box with \emph{open boundaries}
was considered (see e.g.\ den Hollander, Olivieri, and Scoppola~\cite{dHOS00}, den
Hollander, Nardi, Olivieri, and Scoppola~\cite{dHNOS03}, and Bovier, den Hollander,
and Nardi~\cite{BdHN06}).

\medskip
The dynamics of the model will be the continuous-time Markov chain, $(\s_t)_{t \geq 0}$,
with state space $\cX_\b^{(n_\b)}$ whose transition rates are
\be{rate*}
c_\b(\s,\s') =
\left\{\begin{array}{ll}
e^{-\b [H_\b(\s')-H_\b(\s)]_+},
&\mbox{for } \s'=\s^{x,y} \mbox{ for some } x,y \in \L_\b
\mbox{ with } x \sim y,\\
0,
&\mbox{otherwise},
\end{array}
\right.
 \ee
where $\s^{x,y}$ is the configuration obtained from $\s$ by interchanging the values at
sites $x$ and $y$. We refer to this Markov process as \emph{Kawasaki dynamics}. It is
ergodic and reversible with respect to the \emph{canonical} Gibbs measure
\be{mucan}
\mu_\b(\s) = \frac{1}{Z_\b^{(n_\b)}}\,e^{-\b H_\b(\s)}, \qquad \s\in\cX_\b^{(n_\b)},
\ee
where $Z_\b^{(n_\b)}$ is the normalizing partition function. Note that the dynamics preserves
particles, i.e., it is \emph{conservative}.

\begin{figure}[htbp]
\vspace{0.5cm}
\centering
\includegraphics[height=4cm]{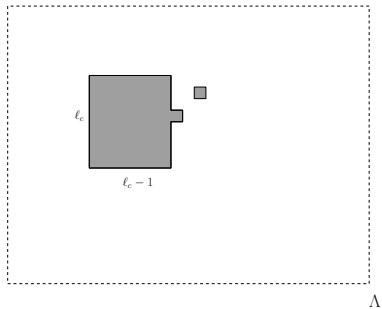}
\caption{\small
A critical droplet for Kawasaki dynamics on $\L$ (= a proto-critical droplet
plus a free particle). The shaded area represents the particles, the non-shaded area
the vacancies (see (\ref{lcdef*})). Note that the shape of the proto-critical droplet
for Kawasaki dynamics is the same as that of the critical droplet for Glauber dynamics.
The proto-critical droplet for Kawasaki dynamics becomes critical when a free particle
is added.}
\label{fig:critdropK}
\end{figure}

\begin{figure}[htbp]
\vspace{0.5cm}
\centering
\includegraphics[height=6cm]{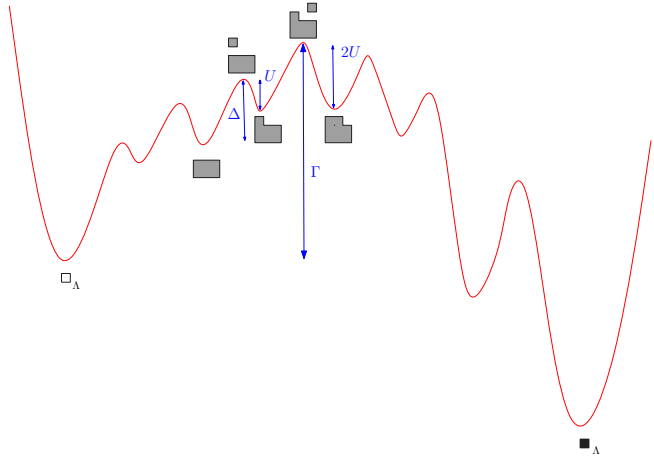}
\caption{\small
A nucleation path from $\square_\L$ to $\blacksquare_\L$ for Kawasaki dynamics
on $\L$ with open boundary. $\G$ in (\ref{Gammadef*}) is the minimal energy barrier
the path has to overcome under the local variant of the grand-canonical Hamiltonian
in (\ref{Hamgc}).}
\label{fig:nuclpath*}
\end{figure}

Kawasaki dynamics exhibits metastable behavior in the regime
\be{metreg*}
U < \Delta < 2U, \qquad \b \to \infty.
\ee
This is again inferred from the behavior of the model in a finite $\b$-independent box
$\L\subset\Z^2$. Let $\square_{\L}$ and $\blacksquare_{\L}$ denote the configurations
where all the sites in $\L$ are vacant, respectively, occupied. For Kawasaki dynamics
on $\L$ with an \emph{open boundary}, where particles are annihilated at rate $1$
and created at rate $e^{-\Delta\b}$, it was shown in den Hollander, Olivieri, and
Scoppola~\cite{dHOS00} and in Bovier, den Hollander, and Nardi~\cite{BdHN06} that,
subject to (\ref{metreg*}) and for the Hamiltonian in (\ref{Hamgc}), the critical
droplets for the crossover from $\square_{\L}$ to $\blacksquare_{\L}$ are the set of
all those configurations where the particles form
\begin{itemize}
\item[(1)]
either an $(\ell_c-2)\times(\ell_c-2)$ square with four bars attached to the
four sides with total length $3\ell_c-3$,
\item[(2)]
or an $(\ell_c-1)\times(\ell_c-3)$ rectangle with four bars attached to the
four sides with total length $3\ell_c-2$,
\end{itemize}
plus a \emph{free particle} anywhere in the box, where
\be{lcdef*}
\ell_c = \left\lceil\frac{U}{2U-\Delta}\right\rceil
\ee
(see Figs.~\ref{fig:critdropK} and \ref{fig:nuclpath*}; for non-degeneracy reasons it
is assumed that $U/(2U-\Delta)\notin\N$).

Let us now return to our setting with finite $\b$-dependent volumes. We define
a reference distance, $L_\b$, as
\be{Lbetadef}
L_\b^2 = e^{(\Delta-\d_\b)\b} = \frac{1}{\rho_\b}\,e^{-\d_\b\b}
\ee
with
\be{dbprop}
\lim_{\beta\to\infty} \d_\b = 0, \qquad \lim_{\b\to\infty} \b\d_\b = \infty,
\ee
i.e., $L_\b$ is marginally below the typical interparticle distance. We assume
$L_\b$ to be odd, and write $B_{L_\b,L_\b}(x)$, $x\in\L_\b$, for the square box
with side length $L_\b$ whose center is $x$.

\bd{Rdef*}
(a) $\cS = \{\s\in\cX_\b^{(n_\b)}\colon\,|\textnormal{supp}[\s] \cap
B_{L_\b,L_\b}(x)|\leq \ell_c(\ell_c-1)+1\,\,\forall\,x\in\L_\b\}$.\\
(b) $\cP=\{\s\in\cS\colon c_\b(\s,\s')>0 \mbox{ for some } \s'\in
\cS^c\}$.\\
(c) $\cC=\{\s'\in\cS^c\colon c_\b(\s,\s')>0 \mbox{ for some } \s\in
\cS\}$.\\
(d) $\cC^-=\{\s\in\cC\colon\,\exists\,x \in \L_\b$ such that $B_{L_\b,L_\b}(x)$
contains a proto-critical droplet plus a free particle at distance $L_\b\}$.\\
(e) $\cC^+$ = the set of configurations obtained from $\cC^-$ by moving the
free particle to a site at distance $2$ from the proto-critical droplet.
\ed

\noindent
As before, we refer to $\cS$, $\cP$ and $\cC$ as the set of \emph{subcritical},
\emph{proto-critical}, respectively, \emph{critical} configurations. Note that,
for every $\s\in\cS$, the number of particles in a box of size $L_\b$ does not
exceed the number of particles in a proto-critical droplet. These particles do
not have to form a cluster or to be near to each other, because the Kawasaki
dynamics brings them together in a time of order $L_\b^2=o(1/\rho_\b)$.

The initial law will again be concentrated on sets $\cS_L\subset\cS$, this time
defined by
\be{SLdefKa}
\cS_L = \big\{\s\in\cX_\b^{(n_\b)}\colon\,|\textnormal{supp}[\s] \cap
B_{L_\b,L_\b}(x)|\leq L\,\,\forall\,x\in\L_\b\big\},
\ee
and $L$ any integer satisfying
\be{L*defKa}
L^* \leq L \leq \ell_c(\ell_c-1)+1 \quad \mbox{ with } \quad
L^* = \min \left\{1 \leq L \leq \ell_c(\ell_c-1)+1\colon\,
\lim_{\b\to\infty} \frac{\mu_\b(\cS_L)}{\mu_\b(\cS)} = 1\right\}.
\ee
In words, $\cS_L$ is the subset of those subcritical configurations for which no
box of size $L_\b$ carries more than $L$ particles, with $L$ again chosen such
that $\cS_L$ is typical within $\cS$ under the Gibbs measure $\mu_\b$ as $\b\to\infty$.

\medskip\noindent
{\bf Remark:}
Note that $\cS_{\ell_c(\ell_c-1)+1}=\cS$. As for Glauber, the value of $L^*$ depends
on how fast $\L_\b$ grows with $\b$. In Appendix \ref{appC.2} we will show that, for
every $1 \leq L \leq \ell_c(\ell_c-1)$, $\lim_{\b\to\infty}\mu_\b(\cS_L)/\mu_\b(\cS)=1$
if and only if $\lim_{\b\to\infty}|\L_\b|e^{-\b(\G_{L+1}-\Delta)}=0$ with $\G_{L+1}$ the
energy needed to create a droplet of $L+1$ particles, closest in shape to a square or
quasi-square, in $B_{L_\b,L_\b}(0)$ under the grand-canonical Hamiltonian on this box.
Thus, if $|\L_\b|=e^{\theta\beta}$, then $L^*= L^*(\theta)=[\ell_c(\ell_c-1)+1] \wedge
\min\{L\in\N\colon\,\Gamma_{L+1}-\Delta>\theta\}$, which increases stepwise from $1$ to
$\ell_c(\ell_c-1)+1$ as $\theta$ increases from $\Delta$ to $\Gamma$ defined in
(\ref{Gammadef*}).

\medskip
Set
\be{Gammadef*}
\Gamma = -U[(\ell_c-1)^2+\ell_c(\ell_c-1)+1]+\Delta[\ell_c(\ell_c-1)+2],
\ee
which is the energy of a critical droplet at a given location with respect to the
\emph{grand-canonical} Hamiltonian given by (\ref{Hamgc}) (see Figs.~\ref{fig:critdropK}
and \ref{fig:nuclpath*}). Put $N = \frac13 \ell_c^2(\ell_c^2-1)$. For $M\in\N$ with $M
\geq \ell_c$, define
\be{DMdef*}
\cD_M = \big\{\s\in\cX_\b\colon\,\exists\,x\in\L_\b \mbox{ such that }
\textnormal{supp}[(\s)] \supset R_{M,M}(x)\big\},
\ee
i.e., the set of configurations containing a supercritical droplet of size $M$. For our
results below to be valid we need to assume that
\be{LambdacondsKa}
\lim_{\b\to\infty} |\L_\b|\,\rho_\b=\infty, \qquad
\lim_{\b\to\infty} |\L_\b|\,e^{-\b\G} = 0.
\ee
This first condition says that the number of particles tends to infinity, and ensures
that the formation of a critical droplet somewhere does not globally deplete the
surrounding gas.

\bt{HNKawasaki}
In the regime {\rm (\ref{metreg*})}, subject to {\rm (\ref{L*defKa})} and {\rm
(\ref{LambdacondsKa})}, the following hold:\\
(a)
\be{crossbdsKawasaki2}
\lim_{\b\to\infty} |\L_\b|\,
\frac{4\pi}{\b\Delta}\,e^{-\b\Gamma}
\,\E_{\nu_{\cS_L}^{(\cS^c\backslash\tilde{\cC})\cup\cC^+}}
\big(\tau_{(\cS^c\backslash\tilde{\cC})\cup\cC^+}\big) = \frac{1}{N}.
\ee
(b)
\be{crossbdsKawasaki3}
\lim_{\b\to\infty} |\L_\b|\,\frac{4\pi}{\b\Delta}\,
e^{-\b\Gamma}
\,\E_{\nu_{\cS_L}^{\cD_M}}\left(\tau_{\cD_M}\right) = \frac{1}{N},
\qquad \forall\,\ell_c \leq M \leq 2\ell_c-1.
\ee
\et

\noindent
The proof of Theorem~\ref{HNKawasaki}, which is the analog of Theorem~\ref{HNGlauber},
will be given in Section~\ref{S4}. Part (a) says that the average time to create a
critical droplet is $[1+o(1)](\b\Delta/4\pi)e^{\b\Gamma}N|\L_\b|$. The factor $\b\Delta/4\pi$
comes from the simple random walk that is performed by the free particle ``from the gas to
the proto-critical droplet'' (i.e., the dynamics goes from $\cC^-$ to $\cC^+$), which slows
down the nucleation. The factor $N$ counts the number of shapes of the proto-critical droplet
(see Bovier, den Hollander, and Nardi~\cite{BdHN06}). Part (b) says that, once the critical
droplet is created, it rapidly grows to a droplet that has twice the size.

\medskip\noindent
{\bf Remarks:}
(1) As for Theorem~\ref{HNGlauber}(c), we expect  Theorem~\ref{HNKawasaki}(b) to hold for
values of $M$ that grow with $\b$ as $M=e^{o(\b)}$. See Section~\ref{S4.2} for more details.\\
(2) In Appendix~\ref{appD} we will show that the \emph{average} probability under
the Gibbs measure $\mu_\b$ of destroying a supercritical droplet and returning to
a configuration in $\cS_L$ is exponentially small in $\b$. Hence, the crossover from
$\cS_L$ to $\cS^c\backslash\tilde{\cC}\cup\cC^+$ represents the true threshold for
nucleation, and Theorem~\ref{HNKawasaki}(a) represents the true nucleation time.\\
(3) It was shown in Bovier, den Hollander, and Nardi~\cite{BdHN06} that the average
crossover time in a finite box $\L$ equals
\be{Kasymp*}
Ke^{\b\Gamma}[1+o(1)], \mbox{ with } K = K(\L,\ell_c)
\sim \frac{\log|\L|}{4\pi}\,\frac{1}{N|\L|},\,\,\L \to \Z^2.
\ee
This matches the $|\L_\b|$-dependence in Theorem~\ref{HNKawasaki}, with the
logarithmic factor in (\ref{Kasymp*}) accounting for the extra factor $\b\Delta$
in Theorem~\ref{HNKawasaki} compared to Theorem~\ref{HNGlauber}. Note that
this factor is particularly interesting, since it says that the \emph{effective
box size} responsible for the formation of a critical droplet is $L_\b$.

\subsection{Outline}
\label{S1.3}
The remainder of this paper is organized as follows. In Section~\ref{S2} we present
a brief sketch of the basic ingredients of the potential-theoretic approach to
metastability. In particular, we exhibit a \emph{relation} between average crossover
times and capacities, and we state \emph{two variational representations} for
capacities, the first of which is suitable for deriving upper bounds and the second
for deriving lower bounds. Section~\ref{S3} contains the proof of our results for the
case of Glauber dynamics. This will be technically relatively easy, and will give a
first flavor of how our method works. In Section~\ref{S4} we deal with Kawasaki dynamics.
Here we will encounter several rather more difficult issues, all coming from the fact
that Kawasaki dynamics is \emph{conservative}. The first is to understand why the
constant $\G$, representing the local energetic cost to create a critical droplet,
involves the grand-canonical Hamiltonian, even though we are working in the canonical
ensemble. This mystery will, of course, be resolved by the observation that the formation
of a critical droplet reduces the entropy of the system: the precise computation of this
entropy loss yields $\Gamma$ via \emph{equivalence of ensembles}. The second problem
is to control the probability of a particle moving from the gas to the proto-critical
droplet at the last stage of the nucleation. This non-locality issue will be dealt with
via upper and lower estimates. Appendices~\ref{appA}--\ref{appD} collect some technical
lemmas that are needed in Sections~\ref{S3}--\ref{S4}.

The extension of our results to higher dimensions is limited only by the combinatorial
problems involved in the computation of the number of critical droplets (which is hard
in the case of Kawasaki dynamics) and of the probability for simple random walk to hit
a critical droplet of a given shape when coming from far. We will not pursue this
generalization here. The relevant results on a $\b$-independent box in $\Z^3$ can
be found in Ben Arous and Cerf~\cite{BC96} (Glauber) and den Hollander, Nardi, Olivieri,
and Scoppola~\cite{dHNOS03} (Kawasaki). For recent overviews on droplet growth in
metastability, we refer the reader to den Hollander~\cite{dH04,dH07} and
Bovier~\cite{B06,B07}. A general overview on metastability is given in the
monograph by Olivieri and Vares~\cite{OV04}.


\section{Basic ingredients of the potential-theoretic approach}
\label{S2}
\setcounter{equation}{0}

The proof of Theorems \ref{HNGlauber} and \ref{HNKawasaki} uses the
\emph{potential-theoretic approach} to metastability developed in
Bovier, Eckhoff, Gayrard and Klein \cite{BEGK02}. This approach is based
on the following three observations. First, most quantities of physical
interest can be represented in term of Dirichlet problems associated with
the generator of the dynamics. Second, the Green function of the dynamics
can be expressed in terms of capacities and equilibrium potentials. Third,
capacities satisfy variational principles that allow for obtaining upper
and lower bounds in a flexible way. We will see that in the current setting
the implementation of these observations provides very sharp results.

\subsection{Equilibrium potential and capacity}
\label{Sepotcap}

The fundamental quantity in the theory is the \emph{equilibrium potential},
$h_{\cA,\cB}$, associated with two non-empty disjoint sets of configurations,
$\cA,\cB\subset\cX$ (= $\cX_\b$ or $\cX_\b^{(n_\b)}$), which probabilistically
is given by
\be{hminprop}
h_{\cA,\cB}(\s) = \left\{\begin{array}{ll}
\P_\s(\tau_\cA<\tau_\cB),
&\mbox{for } \s \in (\cA \cup \cB)^c,\\
1,
&\mbox{for } \s \in \cA,\\
0,
&\mbox{for } \s \in \cB,
\end{array}
\right.
\ee
where
\be{time}
\tau_\cA = \inf \{t>0\colon\, \s_t \in \cA, \s_{t^-}\notin\cA\},
\ee
$(\s_t)_{t \geq 0}$ is the continuous-time Markov chain with state space $\cX$,
and $\P_\s$ is its law starting from $\s$. This function is \emph{harmonic} and
is the unique solution of the \emph{Dirichlet problem}
\be{harmL}
\begin{array}{llll}
(Lh_{\cA,\cB})(\s) &=& 0,  &\s \in (\cA \cup \cB)^c,\\
h_{\cA,\cB}(\s)    &=& 1,  &\s \in \cA,\\
h_{\cA,\cB}(\s)    &=& 0,  &\s \in \cB,
\end{array}
\ee
where the generator is the matrix with entries
\be{gker}
L(\s,\s')= c_\b(\s,\s') - \delta_{\s,\s'}\,c_\b(\s),
\qquad \s,\s'\in\cX,
\ee
with $c_\b(\s)$ the total rate at which the dynamics leaves $\s$,
\be{csigmadef}
c_\b(\s) = \sum_{\s'\in\cX\backslash\{\s\}} c_\b(\s,\s'), \qquad \s\in\cX.
\ee
A related quantity is the \emph{equilibrium measure} on $\cA$, which is defined as
\be{eqm.1}
e_{\cA,\cB}(\s) = -(Lh_{\cA,\cB})(\s), \qquad\s\in\cA.
\ee
The equilibrium measure also has a probabilistic meaning, namely,
\be{eqm.2}
\P_\s(\t_{\cB}<\t_{\cA})= \frac{e_{\cA,\cB}(\s)}{c_\b(\s)}, \qquad\s\in\cA.
\ee
The key object we will work with is the \emph{capacity}, which is defined as
\be{eqm.3}
\CAPA(\cA,\cB) = \sum_{\s\in \cA} \mu_\b(\s)e_{\cA,\cB}(\s).
\ee

\subsection{Relation between crossover time and capacity}
\label{S2.1}

The first important ingredient of the potential-theoretic approach to metastability
is a formula for the average crossover time from $\cA$ to $\cB$. To state this formula,
we define the probability measure $\nu_{\cA}^{\cB}$ on $\cA$ we already referred to in
Section \ref{S1}, namely,
\be{lastexitdef}
\nu_\cA^\cB(\s) = \left\{\begin{array}{ll}
\frac{\mu_\b(\s)e_{\cA,\cB}(\s) }{\CAPA(\cA,\cB)}, &\mbox{for } \s \in \cA,\\
0, &\mbox{for } \s \in \cA^c.
\end{array}
\right.
\ee
The following proposition is proved e.g.\ in Bovier \cite{B07}.

\bp{timecaprel}
For any two non-empty disjoint sets $\cA,\cB\subset\cX$,
\be{crosscaprel}
\sum_{\s\in\cA} \nu_\cA^\cB(\s)\,\E_\s(\tau_\cB)
= \frac{1}{\CAPA(\cA,\cB)}\,\sum_{\s\in\cB^c}
\mu_\b(\s)\,h_{\cA,\cB}(\s).
\ee
\ep

\medskip\noindent
{\bf Remarks:}
(1) Due to (\ref{eqm.2}--\ref{eqm.3}), the probability measure $\nu_\cA^\cB(\s)$
can be written as
\be{lastexitdef.1}
\nu_\cA^\cB(\s) =
\frac{\mu_\b(\s)\,c_\b(\s)}{\CAPA(\cA,\cB)}\,
\P_\s(\tau_\cB<\tau_\cA), \qquad \s\in\cA,
\ee
and thus has the flavor of a \emph{last-exit biased distribution}. Proposition
\ref{timecaprel} explains why our main results on average crossover times stated
in Theorem \ref{HNGlauber} and \ref{HNKawasaki} are formulated for this initial
distribution. Note that
\be{potsumbds}
\mu_\b(\cA) \leq \sum_{\s\in\cB^c} \mu_\b(\s)\,h_{\cA,\cB}(\s) \leq \mu_\b(\cB^c).
\ee
We will see that in our setting $\mu_\b(\cB^c\backslash\cA) = o(\mu_\b(\cA))$ as
$\b\to\infty$, so that the sum in the right-hand side of (\ref{crosscaprel}) is
$\sim\mu_\b(\cA)$ and the computation of the crossover time reduces to the estimation
of $\CAPA(\cA,\cB)$.\\
(2)
For a fixed target set $\cB$, the choice of the starting set $\cA$ is free. It is
tempting to choose $\cA=\{\s\}$ for some $\s\in\cX$. This was done for the case
of a finite $\b$-independent box $\L$. However, in our case (and more generally in
cases where the state space is large) such a choice would give intractable numerators
and denominators in the right-hand side of (\ref{crosscaprel}). As a rule, to make
use of the identity in (\ref{crosscaprel}), $\cA$ must be so large that the harmonic
function $h_{\cA,\cB}$ ``does not change abruptly near the boundary of $\cA$'' for the
target set $\cB$ under consideration.

\medskip
As noted above, average crossover times are essentially governed by capacities. The
usefulness of this observation comes from the computability of capacities, as will be
explained next.

\subsection{The Dirichlet principle: A variational principle for upper bounds}
\label{S2.1.2}

The capacity is a boundary quantity, because $e_{\cA,\cB}>0$ only on the boundary
of $\cA$. The analog of Green's identity relates it to a bulk quantity. Indeed, in
terms of the \emph{Dirichlet form} defined by
\be{Diridef}
\cE(h) = \tfrac12 \sum_{\s,\s'\in\cX}
\mu_\b(\s)c_\b(\s,\s')[h(\s)-h(\s')]^2,
\qquad h\colon\,\cX \to [0,1],
\ee
it follows, via (\ref{hminprop}) and (\ref{eqm.2}--\ref{eqm.3}), that
\be{capdef.0}
\CAPA(\cA,\cB) = \cE(h_{\cA,\cB}).
\ee
Elementary variational calculus shows that the capacity satisfies the
\emph{Dirichlet principle}:

\bp{dp}
For any two non-empty disjoint sets $\cA,\cB\subset\cX$,
\be{capdef}
\CAPA(\cA,\cB) = \min_{{h\colon\,\cX\to [0,1]}
\atop {h|_\cA\equiv 1,h|_\cB\equiv 0}} \cE(h).
\ee
\ep
The importance of the Dirichlet principle is that it yields \emph{computable upper
bounds} for capacities by suitable choices of the test function $h$. In metastable
systems, with the proper physical insight it is often possible to guess a reasonable
test function. In our setting this will be seen to be relatively easy.

\subsection{The Berman-Konsowa principle: A variational principle for lower\\ bounds}
\label{S2.2}

We will describe a little-known variational principle for capacities that is originally due
to Berman and Konsowa \cite{BK90}. Our presentation will follow the argument given in
Bianchi, Bovier, and Ioffe~\cite{BBI}.

In the following it will be convenient to think of $\cX$ as the vertex-set of a graph
$(\cX,\cE)$ whose edge-set $\cE$ consists of all pairs $(\s,\s')$, $\s,\s'\in\cX$, for
which $c_\b(\s,\s')>0$.

\bd{unit-flow}
Given two non-empty disjoint sets $\cA,\cB\subset \cX$, a loop-free non-negative unit
flow, $f$, from $\cA$ to $\cB$ is a function $f\colon\,\cE \rightarrow [0,\infty)$ such
that:\\
(a) $(f(e)>0 \Longrightarrow f(-e)=0)$ $\forall\,e\in\cE$.\\
(b) $f$ satisfies Kirchoff's law:
\be{ab.1}
\sum_{\s'\in\cX} f(\s,\s') = \sum_{\s''\in\cX} f(\s'',\s),
\qquad \forall\,\s\in \cX\backslash (\cA\cup \cB).
\ee
(c) $f$ is normalized:
\be{normalization}
\sum_{\s\in \cA} \sum_{\s'\in\cX} f (\s,\s')
= 1 = \sum_{\s''\in\cX} \sum_{\s\in\cB} f(\s'',\s ).
\ee
(d) Any path from $\cA$ to $\cB$ along edges $e$ such that $f(e)>0$ is self-avoiding.\\
The space of all loop-free non-negative unit flows from $\cA$ to $\cB$ is denoted by
$\U_{\cA,\cB}$.
\ed

A natural flow is the \emph{harmonic flow}, which is constructed from the
equilibrium potential $h_{\cA,\cB}$ as
\be{fstarflow}
f_{\cA,\cB} (\s,\s') = \frac1{\CAPA(\cA,\cB)}\, \mu_\b(\s)c_\b(\s,\s' )
\left[h_{\cA,\cB}(\s)- h_{\cA,\cB}(\s')\right]_+, \qquad \s,\s'\in\cX.
\ee
It is easy to verify that $f_{\cA,\cB}$ satisfies (a--d). Indeed, (a) is obvious, (b)
uses the harmonicity of $h_{\cA,\cB}$, (c) follows from (\ref{eqm.1}) and (\ref{eqm.3}),
while (d) comes from the fact that the harmonic flow only moves in directions where
$h_{\cA,\cB}$ decreases.

A loop-free non-negative unit flow $f$ is naturally associated with a probability
measure $\P^f$ on self-avoiding paths, $\g$. To see this, define $F(\s) = \sum_{\s'\in\cX}
f(\s,\s')$, $\s\in\cX\backslash\cB$. Then $\P^f$ is the Markov chain $(\s_n)_{n\in\N_0}$
with initial distribution $\P^f (\s_0) = F(\s_0)1_\cA(\s_0)$, transition probabilities
\be{ab.3}
q^f (\s,\s' )= \frac{f(\s,\s')}{F(\s)}, \qquad \s\in\cX\backslash\cB,
\ee
such that the chain is stopped upon arrival in $\cB$. In terms of this probability measure,
we have the following proposition (see \cite{BBI} for a proof).

\bp{capMarkov}
Let $\cA,\cB\subset\cX$ be two non-empty disjoint sets. Then, with the notation introduced
above,
\be{ab.5}
\CAPA(\cA,\cB)=\sup_{f\in \U_{\cA,\cB}}\E^f\left(\left[
\sum_{e\in\g}\frac{f(e_l,e_r)}{\mu_\b(e_l)c_\b(e_l,e_r)}\right]^{-1}\right),
\ee
where $e=(e_l,e_r)$ and the expectation is with respect to $\gamma$. Moreover, the supremum
is realized for the harmonic flow $f_{\cA,\cB}$.
\ep

The nice feature of this variational principle is that any flow gives a \emph{computable lower
bound}. In this sense (\ref{capdef}) and (\ref{ab.5}) complement each other. Moreover, since
the harmonic flow is optimal, a good approximation of the harmonic function $h_{\cA,\cB}$ by
a test function $h$ leads to a good approximation of the harmonic flow $f_{\cA,\cB}$ by a test
flow $f$ after putting $h$ instead of $h_{\cA,\cB}$ in (\ref{fstarflow}). Again, in metastable
systems, with the proper physical insight it is often possible to guess a reasonable flow. We
will see in Sections \ref{S3}--\ref{S4} how this is put to work in our setting.


\section{Proof of Theorem \ref{HNGlauber}}
\label{S3}

\setcounter{equation}{0}

\subsection{Proof of Theorem \ref{HNGlauber}(a)}
\label{S3.1}

To estimate the average crossover time from $\cS_L\subset\cS$ to $\cS^c$, we will use
Proposition \ref{timecaprel}. With $\cA=\cS_L$ and $\cB=\cS^c$, (\ref{crosscaprel})
reads
\be{crosscprelSRc}
\sum_{\s\in\cS_L} \nu_{\cS_L}^{\cS^c}(\s)\,\E_\s(\tau_{\cS^c})
= \frac{1}{\CAPA(\cS_L,\cS^c)}\,\sum_{\s\in\cS} \mu_\b(\s)\,h_{\cS_L,\cS^c}(\s).
\ee
The left-hand side is the quantity of interest in (\ref{crossbdsGlauber}). In Sections
\ref{S3.1.1}--\ref{S3.1.2} we estimate $\sum_{\s\in\cS}\mu_\b(\s)h_{\cS_L,\cS^c}(\s)$
and $\CAPA(\cS_L,\cS^c)$. The estimates will show that
\be{rhstot}
\textnormal{r.h.s.}\,{\rm (\ref{crosscprelSRc})}
= \frac{1}{N_1|\L_\b|}\,e^{\b\Gamma}\,[1+o(1)], \qquad \b\to\infty.
\ee

\subsubsection{Estimate of $\sum_{\s\in\cS} \mu_\b(\s)h_{\cS_L,\cS^c}(\s)$}
\label{S3.1.1}

\bl{rhsasymp}
$\sum_{\s\in\cS} \mu_\b(\s) h_{\cS_L,\cS^c}(\s) = \mu_\b(\cS)[1+o(1)]$
as $\b\to\infty$.
\el

\bpr
Write, using (\ref{hminprop}),
\be{rhssumest}
\begin{aligned}
\sum_{\s\in\cS} \mu_\b(\s)h_{\cS_L,\cS^c}(\s)
&= \sum_{\s\in\cS_L} \mu_\b(\s) h_{\cS_L,\cS^c}(\s)
+ \sum_{\s\in\cS\backslash\cS_L} \mu_\b(\s) h_{\cS_L,\cS^c}(\s)\\
&= \mu_\b(\cS_L)
+ \sum_{\s\in\cS\backslash\cS_L} \mu_\b(\s)\P_\s(\tau_{\cS_L}<\tau_{\cS^c}).
\end{aligned}
\ee
The last sum is bounded above by $\mu_\b(\cS\backslash\cS_L)$. But $\mu_\b(\cS\backslash\cS_L)
= o(\mu_\b(\cS))$ as $\b\to\infty$ by our choice of $L$ in (\ref{L*def}).
\epr

\subsubsection{Estimate of $\CAPA(\cS_L,\cS^c)$}
\label{S3.1.2}

\bl{CAPSRcest}
$\CAPA(\cS_L,\cS^c) = N_1\,|\L_\b| e^{-\b\Gamma}\mu_\b(\cS)[1+o(1)]$ as
$\b\to\infty$ with $N_1=4\ell_c$.
\el

\bpr
The proof proceeds via upper and lower bounds.

\medskip\noindent
\underline{Upper bound}:
We use the Dirichlet principle and a test function that is equal to
$1$ on $\cS$ to get the upper bound
\be{Gub1}
\CAPA(\cS_L,\cS^c) \leq \CAPA(\cS,\cS^c)
= \sum_{{\s\in\cS,\s'\in\cS^c} \atop {c_\b(\s,\s')>0}}
\mu_\b(\s)c_\b(\s,\s') = \sum_{{\s\in\cS,\s'\in\cS^c} \atop {c_\b(\s,\s')>0}}
[\mu_\b(\s) \wedge \mu_\b(\s')] \leq \,\mu_\b(\cC),
\ee
where the second equality uses (\ref{rev}) in combination with the fact that
$c_\b(\s,\s') \vee c_\b(\s',\s) = 1$ by (\ref{rate}). Thus, it suffices
to show that
\be{upp2}
\mu_\b(\cC) \leq N_1\,|\L_\b|\,e^{-\b\Gamma}\,[1+o(1)] \qquad \mbox{ as }
\b\to\infty.
\ee
For every $\s\in\cP$ there are one or more rectangles $R_{\ell_c-1,\ell_c}(x)$,
$x=x(\s)\in\cX_\b$, that are filled by $(+1)$-spins in $C_B(\s)$. If $\s' \in\cC$
is such that $\s'=\s^{y}$ for some $y\in\L_\b$, then $\s'$ has a $(+1)$-spin at $y$
situated on the boundary of one of these rectangles. Let
\be{sigmau}
\begin{aligned}
\hat{\cS}(x) &= \big\{\s\in\cS\colon\,
\textnormal{supp}[\s] \subseteq R_{\ell_c-1,\ell_c}(x)\big\},\\
\check{\cS}(x) &=\big\{\s\in\cS\colon\,
\textnormal{supp}[\s]\subseteq [R_{\ell_c+1,\ell_c+2}(x-(1,1))]^c\big\}.
\end{aligned}
\ee

\begin{figure}[htbp]
\vspace{0.5cm}
\centering
\includegraphics[height=3.5cm]{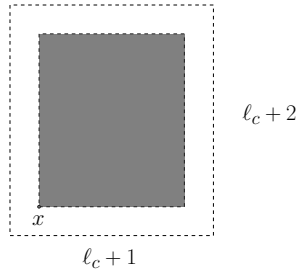}
\caption{\small
$R_{\ell_c-1,\ell_c}(x)$ (shaded box) and $[R_{\ell_c+1,\ell_c+2}(x-(1,1))]^c$
(complement of dotted box).}
\label{fig:rell}
\end{figure}

\medskip\noindent
For every $\s\in\cP$, we have $\s=\hat{\s}\vee\check{\s}$ for some
$\hat{\s}\in\hat{\cS}(x)$ and $\check{\s}\in\check{\cS}(x)$, uniquely
decomposing the configuration into two non-interacting parts inside $R_{\ell_c-1,
\ell_c}(x)$ and $[R_{\ell_c+1,\ell_c+2}(x-(1,1))]^c$ (see Fig.~\ref{fig:rell}). We
have
\be{ener}
H_\b(\s) - H_\b(\boxminus) = [H_\b(\hat{\s}) - H_\b(\boxminus)]
+ [H_\b(\check{\s})-H_\b(\boxminus)].
\ee
Moreover, for any $y\notin\textnormal{supp}[C_B(\s)]$, we have
\be{eflip}
H_\b(\s^{y})\geq H_\b(\s) + 2J - h.
\ee
Hence
\be{muparRcets}
\begin{aligned}
\mu_\b(\cC) &= \frac{1}{Z_\b}\,\sum_{\s\in\cP}\,
\sum_{ {x\in\L_\b} \atop {\s^{x}\in\cC} }\,
e^{-\b H_\b(\s^{x})}\\
&\leq \frac{1}{Z_\b}\,N_1\,e^{-\b[2J-h-H_\b(\boxminus)]}\,
\sum_{x\in\L_\b}\,\sum_{\check{\s}\in\check{\cS}(x)}\,
e^{-\b H_\b(\check{\s})}
\sum_{ {\hat{\s}\in\hat{\cS}(x)} \atop {\hat{\s}\vee\check{\s}\in\cP}}\,
e^{-\b H_\b(\hat{\s})}\\
&\leq [1+o(1)]\,\frac{1}{Z_\b}\, N_1\,|\L_\b|\,e^{-\b\Gamma}\,
\sum_{\check{\s}\in\check{\cS}(0)}\,
e^{-\b H_\b(\check{\s})}\\
&= [1+o(1)]\,N_1\,|\L_\b|\,e^{-\b\Gamma}\,\mu_\b(\check{\cS}(0)),
\end{aligned}
\ee
where the first inequality uses (\ref{ener}--\ref{eflip}), with $N_1 = 2 \times 2\ell_c
= 4\ell_c$ counting the number of critical droplets that can arise from a proto-critical
droplet via a spin flip (see Fig.~\ref{fig:critdropG}), and the second inequality uses
that
\be{Hhatsig}
\hat{\s}\in\hat{\cS}(0),\,\hat{\s}\vee\check{\s}\in\cP
\Longrightarrow H_\b(\hat{\s}) \geq H_\b(R_{\ell_c-1,\ell_c}(0))
= \Gamma-(2J-h)+H_\b(\boxminus)
\ee
with equality in the right-hand side if and only if $\textnormal{supp}[\hat{\s}] =
R_{\ell_c-1,\ell_c}(0)$. Combining (\ref{Gub1}) and (\ref{muparRcets}) with the
inclusion $\check{\cS}(0)\subset\cS$, we get the upper bound in (\ref{upp2}).

\begin{figure}[htbp]
\vspace{0.5cm}
\centering
\includegraphics[height=1.5cm]{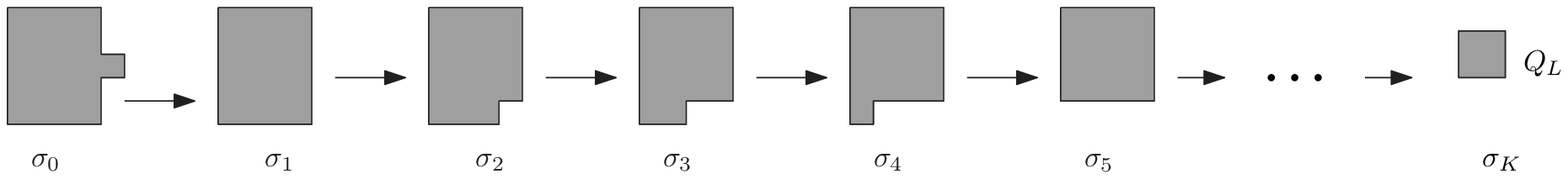}
\caption{\small
Canonical order to break down a critical droplet. }
\label{fig:droptear}
\end{figure}

\medskip\noindent
\underline{Lower bound}:
We exploit Proposition \ref{capMarkov} by making a judicious choice
for the flow $f$. In fact, in the Glauber case this choice will be
simple: with each configuration $\s\in\cS_L$ we associate a configuration
in $\cC\subset\cS^c$ with a unique critical droplet and a flow that, from each
such configuration, follows a unique \emph{deterministic path} along which this
droplet is broken down in the \emph{canonical order} (see Fig.~\ref{fig:droptear})
until the set $\cS_L$ is reached, i.e., a square or quasi-square droplet with
label $L$ is left over (recall (\ref{qsqdef}--\ref{SLdef})).

\begin{figure}[htbp]
\vspace{0.5cm}
\centering
\includegraphics[height=2cm]{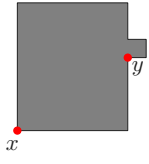}
\caption{\small
The critical droplet $P_{(y)}(x)$.}
\label{fig:Pyx}
\end{figure}

\begin{figure}[htbp]
\vspace{0.5cm}
\centering
\includegraphics[height=6cm]{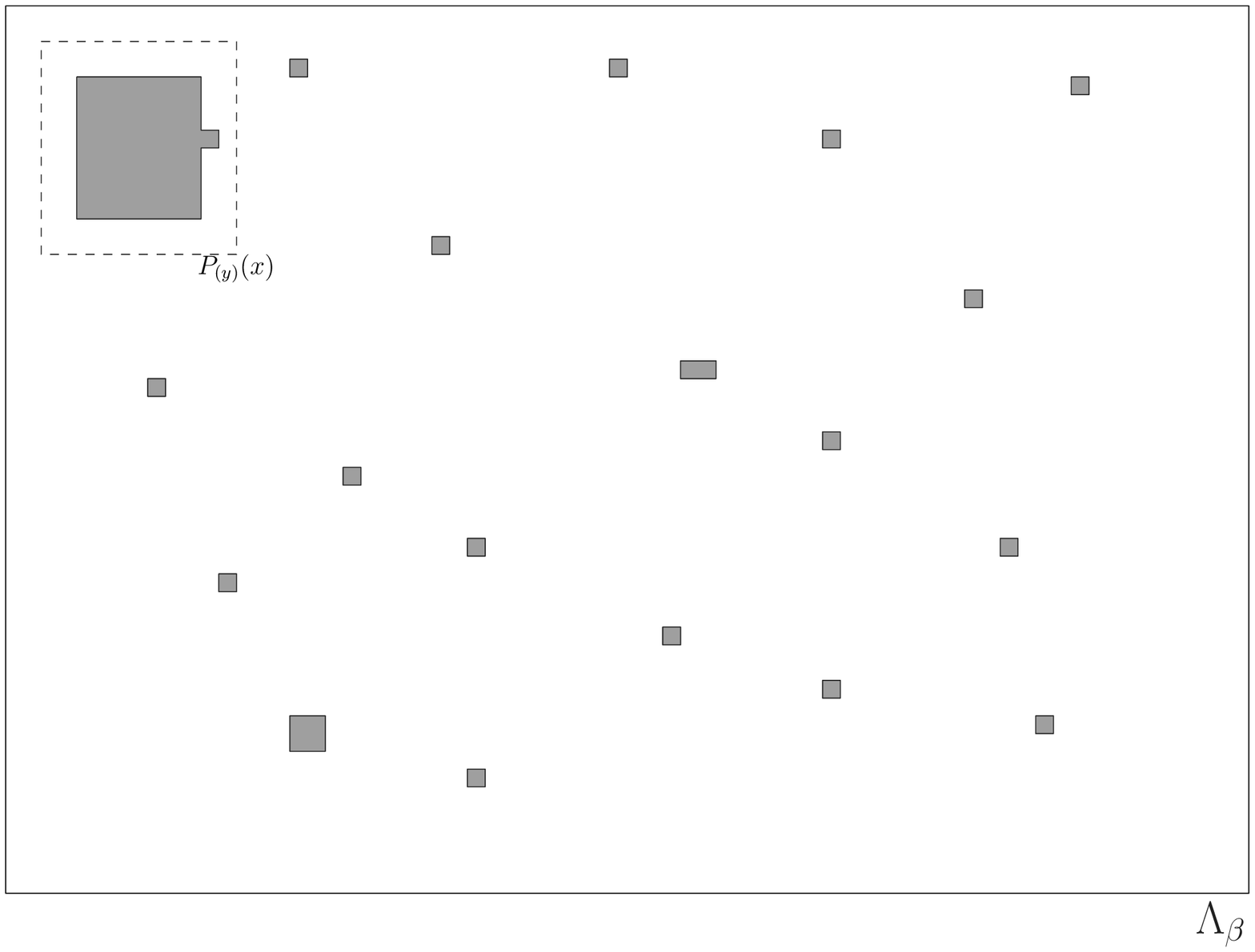}
\caption{\small
Going from $\cS_L$ to $\cC_L$ by adding a critical droplet
$P_{(y)}(x)$ somewhere in $\L_\b$.}
\label{fig:insertcritdrop}
\end{figure}

Let $f(\b)$ be such that
\be{fprop}
\lim_{\b\to\infty} f(\b) = \infty, \quad
\lim_{\b\to\infty} \frac{1}{\b} \log f(\b) = 0, \quad
\lim_{\b\to\infty} |\L_\b|/f(\b) = \infty,
\ee
and define
\be{Wdef}
\cW = \big\{\s\in\cS\colon\,|\textnormal{supp}[\s]|\leq |\L_\b|/f(\b)\big\}.
\ee
Let $\cC_L\subset\cC\subset\cS^c$ be the set of configurations obtained by picking
any $\s\in\cS_L\cap\cW$ and adding somewhere in $\L_\b$ a critical droplet at
distance $\geq 2$ from $\textnormal{supp}[\s]$. Note that the density restriction
imposed on $\cW$ guarantees that adding such a droplet is possible almost everywhere
in $\L_\b$ for $\b$ large enough. Denoting by $P_{(y)}(x)$ the critical droplet
obtained by adding a protuberance at $y$ along the longest side of the rectangle
$R_{\ell_c-1,\ell_c}(x)$, we may write
\be{rcap}
\cC_L = \big\{\s \cup P_{(y)}(x)\colon\,\s\in\cS\cap\cW,\,x,y\in\L_\b,\,
(x,y) \bot \s\big\},
\ee
where $(x,y) \bot \s$ stands for the restriction that the critical droplet $P_{(y)}(x)$
is not interacting with ${\rm supp}[\s]$, which implies that $H_\b(\s\cup P_{(y)}(x))
= H_\b(\s)+\Gamma$ (see Figs.~\ref{fig:Pyx} and \ref{fig:insertcritdrop}).

Now, for each $\s\in \cC_L$, we let $\g_\s=(\g_\s(0),\g_\s(1),\dots,\g_\s(K))$
be the canonical path from $\s=\g_\s(0)$ to $\cS_L$ along which the critical droplet
is broken down, where $K=v(2\ell_c-3)-v(L)$ with
\be{vLdef}
v(L) = |Q_L(0)|
\ee
(recall (\ref{qsqdef})). We will choose our flow such that
\be{flow.1}
\begin{aligned}
&f(\s',\s'')\\
&=
\begin{cases}
\nu_0(\s),
&\hbox{\rm if}\,\, \s'=\s,\,\s''=\g_\s(1)\,\hbox{\rm for some}\,\,\s\in \cC_L,\\
\sum_{\tilde\s\in\cC_L} f(\g_{\tilde \s}(k-1),\g_\s(k)),
&\hbox{\rm if}\,\, \s'=\g_\s(k),\,\s''=\g_\s(k+1)\,\hbox{\rm for some}\,\,k\geq 1,\,
\s\in \cC_L,\\
0,
&\hbox{\rm otherwise}.
\end{cases}
\end{aligned}
\ee
Here, $\nu_0$ is some initial distribution on $\cC_L$ that will turn out to be
arbitrary as long as its support is all of $\cC_L$.

We see from (\ref{flow.1}) that the flow increases whenever paths merge. In our case
this happens only after the first step, when the protuberance at $y$ is removed. Therefore
we get the explicit form
\be{flow.1.1}
f(\s',\s'') =
\begin{cases}
\nu_0(\s),
&\hbox{\rm if}\,\s'=\s,\,\s''=\g_\s(1)\,\hbox{\rm for some}\,\,\s\in \cC_L,\\
C\nu_0(\s),
&\hbox{\rm if}\,\s'=\g_\s(k),\,\s''=\g_\s(k+1)\,\hbox{\rm for some}\,\,k\geq 1,\,\s\in \cC_L,\\
0,
&\hbox{\rm otherwise},
\end{cases}
\ee
where $C=2\ell_c$ is the number of possible positions of the protuberance on the
proto-critical droplet (see Fig.~\ref{fig:droptear}). Using Proposition~\ref{capMarkov},
we therefore have
\be{lb2}
\begin{aligned}
\CAPA(\cS_L,\cS^c) &= \CAPA(\cS^c,\cS_L) \geq \CAPA(\cC_L,\cS_L)\\
&\geq \sum_{\s\in\cC_L} \nu_0(\s)
\left[
\sum_{k=0}^{K-1}
\frac{f(\g_\s(k),\g_\s(k+1))}
{\mu_\b(\g_\s(k))c_\b(\g_\s(k),\g_\s(k+1))}
\right]^{-1}\\
&= \sum_{\s\in\cC_L}
\left[\frac{1}{\mu_\b(\s)c_\b(\g_\s(0),\g_\s(1))}
+ \sum_{k=1}^{K-1} \frac{C}{\mu_\b(\g_\s(k))c_\b(\g_\s(k),\g_\s(k+1))}
\right]^{-1}.
\end{aligned}
\ee
Thus, all we have to do is to control the sum between square brackets.

\begin{figure}[htbp]
\centering
\includegraphics[height=5cm]{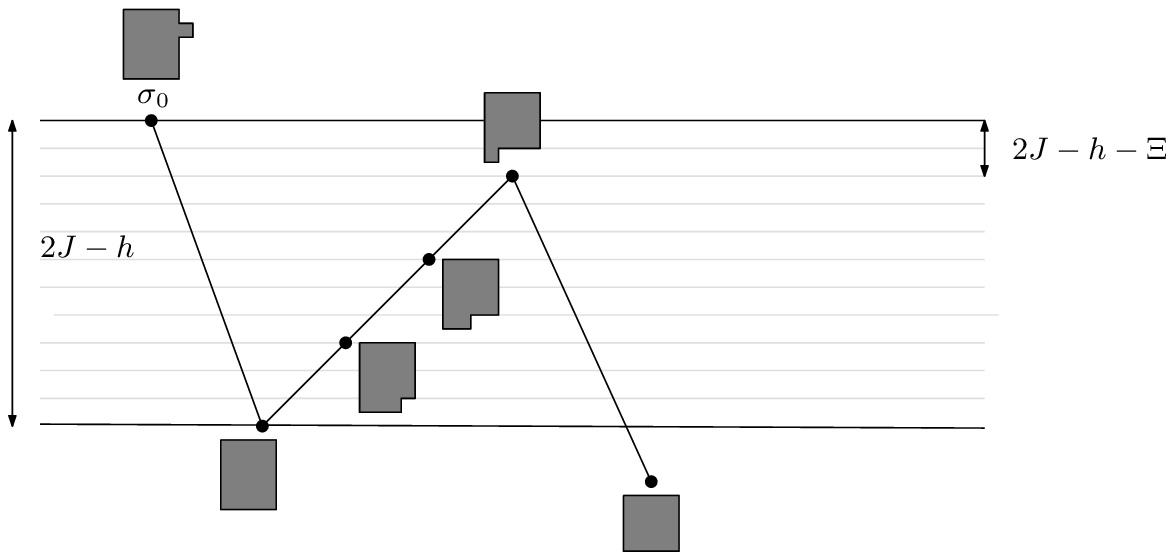}
\caption{\small
Visualization of (\ref{kterms}).}
\label{fig:gap}
\end{figure}

Because $c_\b(\g_\s(0),\g_\s(1))=1$ (removing the protuberance lowers the energy), the
term with $k=0$ equals $1/\mu_\b(\s)$. To show that the terms with $k \geq 1$ are of
higher order, we argue as follows. Abbreviate $\Xi=h(\ell_c-2)$. For every $k\geq 1$
and $\s(0)\in\cC_L$, we have (see Fig.~\ref{fig:gap} and recall (\ref{Gibbs}--\ref{rate}))
\be{kterms}
\mu_\b(\g_\s(k))c_\b(\g_\s(k),\g_\s(k+1))
= \frac{1}{Z_\b}\,e^{-\b[H_\b(\g_\s(k)) \vee H_\b(\g_\s(k+1))]}
\geq \mu_\b(\s_0)\,e^{\b[2J-h-\Xi]}= \mu_\b(\s) e^{\b\delta},
\ee
where $\delta=2J-h-\Xi=2J-h(\ell_c-1)>0$ (recall (\ref{lcdef})). Therefore
\be{sums.2}
\sum_{k=1}^{K-1}
\frac{C}
{\mu_\b(\g_\s(k))c_\b(\g_\s(k),\g_\s(k+1))} \leq  \frac{1}{\mu_\b(\s)}CKe^{-\d\b},
\ee
and so from (\ref{lb2}) we get
\be{CAPlbres}
\CAPA(\cS_L,\cS^c) \geq \sum_{\s\in \cC_L} \frac{\mu_\b(\s)}{1+ CKe^{-\b\delta}}
= \frac{\mu_\b(\cC_L)}{1+ CKe^{-\b\delta}}
= [1+o(1)]\,\mu_\b(\cC_L).
\ee

The last step is to estimate, with the help of (\ref{rcap}),
\be{muRc*est}
\begin{aligned}
\mu_\b(\cC_L)
&= \frac{1}{Z_\b} \sum_{\s\in\cC_L} e^{-\b H_\b(\s)}
= \frac{1}{Z_\b} \sum_{\s\in\cS_L\cap\cW}\,\,
\sum_{ {x,y\in\L_\b} \atop {(x,y) \bot \s} }
e^{-\b H_\b(\s\cup P_{(y)}(x))}\\
&= e^{-\b\Gamma} \frac{1}{Z_\b} \sum_{\s\in\cS_L\cap\cW} e^{-\b H_\b(\s)}
\sum_{ {x,y\in\L_\b} \atop {(x,y) \bot \s} } 1\\
&\geq e^{-\b\Gamma}\,\mu_\b(\cS_L\cap\cW)\,
N_1\,|\L_\b|\,[1-(\ell_c+1)^2/f(\b)].
\end{aligned}
\ee
The last inequality uses that $|\L_\b|(\ell_c+1)^2/f(\b)$ is the maximal number
of sites in $\L_\b$ where it is not possible to insert a non-interacting critical
droplet (recall (\ref{Wdef}) and note that a critical droplet fits inside an
$\ell_c\times\ell_c$ square). According to Lemma~\ref{wses} in Appendix~\ref{appA},
we have
\be{Wasymp}
\mu_\b(\cS_L\cap\cW)=\mu_\b(\cS_L)[1+o(1)],
\ee
while conditions (\ref{SLdef}--\ref{L*def}) imply that $\mu_\b(\cS_L)=\mu_\b(\cS)[1+o(1)]$.
Combining the latter with (\ref{CAPlbres}--\ref{muRc*est}), we obtain the desired lower bound.
\epr

\subsection{Proof of Theorem \ref{HNGlauber}(b)}
\label{S3.2}

We use the same technique as in Section~\ref{S3.1}, which is why we only give a
sketch of the proof.

To estimate the average crossover time from $\cS_L\subset\cS$ to $\cS^c\backslash\cC$,
we will use Proposition \ref{timecaprel}. With $\cA=\cS_L$ and $\cB=\cS^c\backslash\cC$,
(\ref{crosscaprel}) reads
\be{crosscprelSRcalt}
\sum_{\s\in\cS_L} \nu_{\cS_L}^{\cS^c\backslash\cC}(\s)\,
\E_\s(\tau_{\cS^c\backslash\cC})
= \frac{1}{\CAPA(\cS_L,\cS^c\backslash\cC)}\,\sum_{\s\in\cS\cup\cC}
\mu_\b(\s)\,h_{\cS_L,\cS^c\backslash\cC}(\s).
\ee
The left-hand side is the quantity of interest in (\ref{crossbdsGlauber2}).

In Sections \ref{S3.3.1}--\ref{S3.3.2} we estimate $\sum_{\s\in\cS\cup\cC} \mu_\b(\s)
h_{\cS_L,\cS^c\backslash\cC}(\s)$ and $\CAPA(\cS_L,\cS^c\backslash\cC)$. The estimates
will show that
\be{rhstotalt}
\textnormal{r.h.s.}\,{\rm (\ref{crosscprelSRcalt})}
= \frac{1}{N_2|\L_\b|}\,e^{\b\Gamma}\,[1+o(1)], \qquad \b\to\infty.
\ee

\subsubsection{Estimate of $\sum_{\s\in\cS\cup\cC} \mu_\b(\s)
h_{\cS_L,\cS^c\backslash \cC}(\s)$}
\label{S3.3.1}

\bl{rhsasymp*}
$\sum_{\s\in\cS\cup\cC} \mu_\b(\s) h_{\cS_L,\cS^c\backslash \cC}(\s)
= \mu_\b(\cS)[1+o(1)]$ as $\b\to\infty$.
\el

\bpr
Write, using (\ref{hminprop}),
\be{rhssumest2}
\sum_{\s\in\cS\cup\cC} \mu_\b(\s) h_{\cS_L,\cS^c\backslash \cC}(\s)
= \mu_\b(\cS_L) + \sum_{\s\in(\cS\backslash\cS_L)\cup\cC}
\mu_\b(\s)\P_\s(\tau_{\cS_L}<\tau_{\cS^c\backslash\cC}).
\ee
The last sum is bounded above by $\mu_\b(\cS\backslash\cS_L)+\mu_\b(\cC)$. As before,
$\mu_\b(\cS\backslash\cS_L)=o(\mu_\b(\cS))$ as $\b\to\infty$. But (\ref{LambdacondsKa})
and (\ref{muparRcets}) imply that $\mu_\b(\cC)=o(\mu_\b(\cS))$ as $\b\to\infty$.
\epr

\subsubsection{Estimate of $\CAPA(\cS_L,\cS^c\backslash\cC)$}
\label{S3.3.2}

\bl{CAPSRcest2}
$\CAPA(\cS,\cS^c\backslash\cC) = N_2\,|\L_\b| e^{-\b\Gamma}\mu_\b(\cS)[1+o(1)]$ as
$\b\to\infty$ with $N_2=\frac{4}{3}(2\ell_c-1)$.
\el

\bpr
The proof is similar as that of Lemma~\ref{CAPSRcest}, except that it takes care of the
transition probabilities away from the critical droplet.

\medskip\noindent
\underline{Upper bound}:
Recalling (\ref{Diridef}--\ref{capdef}) and noting that Glauber dynamics does not
allow transitions within $\cC$, we have, for all $h\colon\,\cC \to [0,1]$,
\be{gla2}
\CAPA(\cS_L,\cS^c\backslash\cC) \leq \CAPA(\cS,\cS^c\backslash\cC)
\leq \sum_{\s\in\cC}\mu_\b(\s)\big[\hat{c}_{\s}(h(\s)-1)^2
+\check{c}_{\s}(h(\s)-0)^2\big],
\ee
where $\hat{c}_{\s}=\sum_{\eta\in\cS}c_\b(\s,\eta)$ and $\check{c}_{\s}=\sum_{\eta\in
\cS^c\backslash\cC}c_\b(\s,\eta)$. The quadratic form in the right-hand side of $(\ref{gla2})$
achieves its minimum for $h(\s)=\hat{c}_{\s}/(\hat{c}_{\s}+\check{c}_{\s})$, so
\be{gla3}
\CAPA(\cS_L,\cS^c\backslash\cC) \leq \sum_{\s\in\cC} C_\s\,\mu_\b(\s)
\ee
with $C_\s=\hat{c}_{\s}\check{c}_{\s}/(\hat{c}_{\s}+\check{c}_{\s})$.
We have
\be{muparRcets2}
\begin{aligned}
\sum_{\s\in\cC} C_\s\,\mu_\b(\s)
&= \frac{1}{Z_\b}\,\sum_{\s\in\cP}\,
\sum_{ {x\in\L_\b} \atop {\s^{x}\in\cC} }\,C_{\s^x}\,
e^{-\b H_\b(\s^{x})}\\
&= e^{-\b(2J-h)}\,\frac{1}{Z_\b}\,
\sum_{\s\in\cP}\,e^{-\b H_\b(\s)}
\,\,2\left(\tfrac{1}{2}\,4 +\tfrac{2}{3}(2\ell_c-4)\right)\\
&=e^{-\b(2J-h)}\,\mu_\b(\cP)\,N_2 = \frac{1}{N_1}\,\mu_\b(\cC)\,N_2,
\end{aligned}
\ee
where in the second line we use that $C_\s=\frac12$ if $\s$ has a protuberance in a
corner ($2 \times 4$ choices) and $C_\s=\frac23$ otherwise ($2 \times (2\ell_c-4)$
choices).

\begin{figure}[htbp]
\vspace{0.5cm}
\centering
\includegraphics[height=2.5cm]{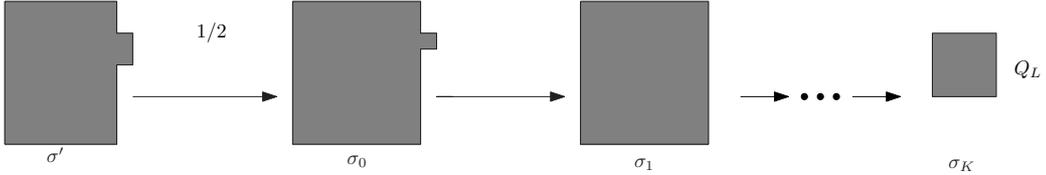}
\caption{\small
Canonical order to break down a proto-critical droplet plus a double protuberance.
In the first step, the double protuberance has probability $\frac12$ to be broken down in either
of the two possible ways. The subsequent steps are deterministic as in Fig.~\ref{fig:droptear}.
}
\label{fig:droptear2}
\end{figure}

\medskip\noindent
\underline{Lower bound}:
In analogy with (\ref{rcap}), denoting by $P_{(y)}^2(x)$ the droplet obtained by adding a
double protuberance at $y$ along the longest side of the rectangle $R_{\ell_c-1,\ell_c}(x)$,
we define the set $\cD_L \subset \cS^c\backslash\cC$ by
\be{rcap2}
\cD_L = \{\s \cup P_{(y)}^2(x)\colon\,\s\in\cS_L\cap\cW,\,
x,y\in\L_\b,\,(x,y) \bot \s\}.
\ee
As in (\ref{flow.1}), we may choose any starting measure on $\cD_L$. We choose the flow
as follows. For the first step we choose
\be{new-flow.1}
f(\s',\s) = \tfrac12\,\nu_0(\s), \qquad \s'\in\cD_L,\,\s\in\cC_L,
\ee
which reduces the double protuberance to a single protuberance (compare (\ref{rcap})
and (\ref{rcap2})). For all subsequent steps we follow the deterministic paths $\g_\s$
used in Section~\ref{S3.1.2}, which start from $\g_\s(0)=\s$. Note, however, that we get
different values for the flows $f(\g_\s(0),\g_\s(1))$ depending on whether the protuberance
sits in a corner or not. In the former case, it has only one possible antecedent, and so
\be{flow-corner.1}
f(\g_\s(0),\g_\s(1)) = \tfrac12\,\nu_0(\s),
\ee
while in the latter case it has two antecedents, and so
\be{flow-corner.2}
f(\g_\s(0),\g_\s(1)) = \nu_0(\s).
\ee
This time the terms $k=0$ and $k=1$ are of the same order while, as in (\ref{sums.2}), all
the subsequent steps give a contribution that is a factor $O(e^{-\d\b})$ smaller. Indeed,
in analogy with (\ref{lb2}) we obtain, writing $\s\sim\s'$ when $c_\b(\s',\s)>0$,
\be{jen}
\begin{aligned}
&\CAPA(\cS_L,\cS^c\backslash\cC) = \CAPA(\cS^c\backslash\cC,\cS_L)
\geq \CAPA(\cD_L,\cS_L)\\
&\geq
\sum_{\s'\in \cD_L} \tfrac12 \sum_{ {\s\in\cC_L} \atop {\s\sim\s'}}
\left[ \frac{f(\s',\s)}{\mu_\b(\s)}+\frac {f(\s,\g_\s(1))}{\mu_\b(\s)}
+\sum_{k=1}^{K-1}\frac{f(\g_\s(k),\g_\s(k+1))}
{\mu_\b(\g_\s(k))c_\b(\g_\s(k),\g_\s(k+1))} \right]^{-1}\\
&\geq \sum_{\s'\in \cD_L} \tfrac12 \sum_{{\s\in\cC_L} \atop {\s\sim\s'}}
\mu_\b(\s) \left[ {f(\s',\s)}+ {f(\s,\g_\s(1))} + CKe^{-\b\d}\right]^{-1}\\
&=[1+o(1)]\,\mu_\b(\cC_L)
\left(
\frac{2\ell_c-4}{2\ell_c}\,\frac{1}{1+\tfrac12}
+\frac 12\,\frac{4}{2\ell_c}\,\frac{1}{\tfrac12+\tfrac12}\right)\\
&= [1+o(1)]\,\mu_\b(\cC_L)\,\frac{N_2}{N_1}.
\end{aligned}
\ee
Using (\ref{muRc*est}) and the remarks following it, we get the desired lower bound.
\epr

\subsection{Proof of Theorem \ref{HNGlauber}(c)}
\label{S3.3}

Write
\be{eqpotMdrop}
\begin{aligned}
\sum_{\s\in\cD_M^c} \mu_\b(\s) h_{\cS_L,\cD_M}(\sigma)
&= \sum_{\s\in\cS_L} \mu_\b(\s) h_{\cS_L,\cD_M}(\sigma)
+ \sum_{\s\in\cD_M^c\backslash\cS_L} \mu_\b(\s) h_{\cS_L,\cD_M}(\sigma)\\
&= \mu_\b(\cS_L) + \sum_{\s\in\cD_M^c\backslash\cS_L} \mu_\b(\s)
\P_\s(\t_{\cS_L} < \t_{\cD_M}).
\end{aligned}
\ee
The last sum is bounded above by $\mu_\b(\cS\backslash\cS_L) + \mu_\b(\cD_M^c\backslash\cS)$.
But $\mu_\b(\cS\backslash\cS_L) = o(\mu_\b(\cS))$ as $\b\to\infty$ by our choice of $L$ in
(\ref{L*def}), while $\mu_\b(\cD_M^c\backslash\cS)= o(\mu_\b(\cS))$ as $\b\to\infty$ because
of the restriction $\ell_c \leq M 2\ell_c-1$. Indeed, under that restriction the energy of
a square droplet of size $M$ is strictly larger than the energy of a critical droplet.

\bpr
The proof of Theorem~\ref{HNGlauber}(c) follows along the same lines as that of
Theorems~\ref{HNGlauber}(a--b) in Sections~\ref{S3.1}--\ref{S3.2}. The main point
is to prove that $\CAPA(\cS_L,\cD_M)=[1+o(1)]\CAPA(\cS_L,\cS^c\backslash\cC)$.
Since $\CAPA(\cS_L,\cD_M)\leq\CAPA(\cS_L,\cS^c\backslash\cC)$, which was estimated in
Section ~\ref{S3.2}, we need only prove a lower bound on $\CAPA(\cS_L,\cD_M)$. This
is done by using a flow that breaks down an $M \times M$ droplet to a square or
quasi-square droplet $Q_L$ in the canonical way, which takes $M^2-v(L)$ steps (recall
Fig.~\ref{fig:droptear} and (\ref{vLdef})). The leading terms are still the proto-critical
droplet with a single and a double protuberance. To each $M\times M$ droplet is associated a
unique critical droplet, so that the pre-factor in the lower bound is the same as in
the proof of Theorem~\ref{HNGlauber}(b).

Note that we can even allow $M$ to grow with $\b$ as $M=e^{o(\b)}$. Indeed,
(\ref{fprop}--\ref{Wdef}) show that there is room enough to add a droplet of size
$e^{o(\b)}$ almost everywhere in $\L_\b$, and the factor $M^2e^{-\d\b}$ replacing
$Ke^{-\d\b}$ in (\ref{CAPlbres}) still is $o(1)$.
\epr


\section{Proof of Theorem \ref{HNKawasaki}}
\label{S4}

\setcounter{equation}{0}

\subsection{Proof of Theorem \ref{HNKawasaki}(a)}
\label{S4.1}

\subsubsection{Estimate of $\sum_{\s\in\cS\cup(\tilde{\cC}\backslash\cC^+)}
\mu_\b(\s) h_{\cS_L,(\cS^c\backslash\tilde{\cC})\cup\cC^+}(\s)$}
\label{S4.1.1}

\bl{rhsasympKa}
$\sum_{\s\in\cS\cup(\tilde{\cC}\backslash\cC^+)}
\mu_\b(\s) h_{\cS_L,(\cS^c\backslash\tilde{\cC})\cup\cC^+}(\s) = \mu_\b(\cS)[1+o(1)]$
as $\b\to\infty$.
\el

\bpr
Write, using (\ref{hminprop}),
\be{rhssumestKa}
\begin{aligned}
&\sum_{\s\in\cS\cup(\tilde{\cC}\backslash\cC^+)} \mu_\b(\s)
h_{\cS_L,(\cS^c\backslash\tilde{\cC})\cup\cC^+}(\s)\\
&\qquad = \mu_\b(\cS_L)
+ \sum_{\s\in(\cS\backslash\cS_L)\cup(\tilde{\cC}\backslash\cC^+)}
\mu_\b(\s) \P_\s\big(\tau_{\cS_L}<\tau_{(\cS^c\backslash\tilde{\cC})\cup\cC^+}\big).
\end{aligned}
\ee
The last sum is bounded above by $\mu_\b(\cS\backslash\cS_L)+\mu_\b(\tilde{\cC}\backslash\cC^+)$.
But $\mu_\b(\cS\backslash\cS_L)= o(\mu_\b(\cS))$ as $\b\to\infty$ by our choice of $L$
in (\ref{L*defKa}). In Lemma~\ref{GibbscC+} in Appendix \ref{appB.3} we will show that
$\mu_\b(\tilde{\cC}\backslash\cC^+) = o(\mu_\b(\cS))$ as $\b\to\infty$.

\epr

\subsubsection{Estimate of $\CAPA(\cS_L,(\cS^c\backslash\tilde{\cC})\cup\cC^+)$}
\label{S4.1.2}

\bl{CAPSRcest*}
$\CAPA(\cS_L,\cS^c\backslash\tilde{\cC})\cup\cC^+)
= N|\L_\b|\,\frac{4\pi}{\b\Delta}\,e^{-\b\Gamma}
\mu_\b(\cS)[1+o(1)]$ as $\b\to\infty$ with $N=\frac13\ell_c^2(\ell_c^2-1)$.
\el

\bpr
The argument is in the same spirit as that in Section~\ref{S3.1.2}. However, a number
of additional hurdles need to be taken that come from the conservative nature of Kawasaki
dynamics. The proof proceeds via upper and lower bounds, and takes up quite a bit of space.

\medskip\noindent
\underline{Upper bound}:
The proof comes in 7 steps.

\begin{figure}[htbp]
\vspace{0.5cm}
\centering
\includegraphics[height=4cm]{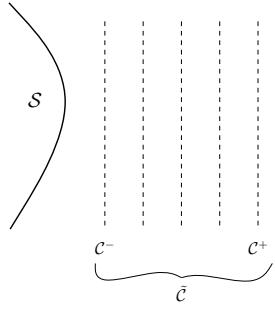}
\caption{\small
Schematic picture of the sets $\cS$, $\cC^-$, $\cC^+$ defined in
Definition~\ref{Rdef*} and the set $\tilde{\cC}$ interpolating between
$\cC^-$ and $\cC^+$.}
\label{fig:SCtilde}
\end{figure}

\medskip\noindent
{\bf 1.\ Proto-critical droplet and free particle.}
Let $\tilde{\cC}$ denote the set of configurations ``interpolating'' between $\cC^-$ and
$\cC^+$, in the sense that the free particle is somewhere between the boundary of the
proto-critical droplet and the boundary of the box of size $L_\b$ around the proto-critical
droplet (see Fig.~\ref{fig:SCtilde}). Then we have
\be{upkaw}
\begin{aligned}
\CAPA(\cS_L,(\cS^c\backslash\tilde{\cC})\cup\cC^+)
&\leq  \CAPA(\cS\cup\cC^-,(\cS^c\backslash\tilde{\cC})\cup\cC^+)\\
&= \min_{ {h\colon\,\cX_\b^{(n_\b)} \to [0,1]} \atop
{h|_{\cS\cup\cC^-} \equiv 1,\,h|_{(\cS^c\backslash\tilde{\cC})\cup\cC^+} \equiv 0} }
\tfrac12 \sum_{\s,\s'\in\cX_\b^{(n_\b)}} \mu_\b(\s)c_\b(\s,\s')\,[h(\s)-h(\s')]^2.
\end{aligned}
\ee
Split the right-hand side into a contribution coming from $\s,\s'\in\tilde{\cC}$ and the rest:
\be{upkawsplit}
\textnormal{r.h.s.} (\ref{upkaw}) = I + \gamma_1(\b),
\ee
where
\be{Idef}
I = \min_{ {h\colon\,\tilde{\cC} \to [0,1]}
\atop {h|_{\cC^-} \equiv 1,\,h|_{\cC^+} \equiv 0} }
\tfrac12 \sum_{\s,\s'\in\tilde{\cC}} \mu_\b(\s)c_\b(\s,\s')\,[h(\s)-h(\s')]^2
\ee
and $\gamma_1(\b)$ is an error term that will be estimated in Step 7. This term will turn out
to be small because $\mu_\b(\s)c_\b(\s,\s')$ is small when either $\s\in \cX_\b^{(n_\b)}\backslash
\tilde{\cC}$ or $\s'\in \cX_\b^{(n_\b)}\backslash\tilde{\cC}$. Next, partition $\tilde{\cC}$, $\cC^-$,
$\cC^+$ into sets $\tilde{\cC}(x)$, $\cC^-(x)$, $\cC^+(x)$, $x\in\L_\b$, by requiring that the
lower-left corner of the proto-critical droplet is in the center of the box $B_{L_\b,L_\b}(x)$.
Then, because $c_\b(\s,\s')=0$ when $\s\in\tilde{\cC}(x)$ and $\s'\in \tilde{\cC}(x')$ for some
$x \neq x'$, we may write
\be{upkawII}
I = |\L_\b|\, \min_{ {h\colon\,\tilde{\cC}(0) \to [0,1]}
\atop {h|_{\cC^-(0)} \equiv 1,\,h|_{\cC^+(0)} \equiv 0} }
\tfrac12 \sum_{\s,\s'\in\tilde{\cC}(0)} \mu_\b(\s)c_\b(\s,\s')\,[h(\s)-h(\s')]^2.
\ee

\medskip\noindent
{\bf 2.\ Decomposition of configurations.}
Define (compare with (\ref{sigmau}))
\be{cchec}
\begin{aligned}
\hat{\cC}(0) &= \big\{\s \1_{B_{L_\b,L_\b}(0)}\colon\,\s\in\tilde{\cC}(0)\big\},\\
\check{\cC}(0) &=\big\{\s \1_{[B_{L_\b,L_\b}(0)]^c}\colon\,\s\in\tilde{\cC}(0)\big\}.
\end{aligned}
\ee
Then every $\s\in\tilde{\cC}(0)$ can be uniquely decomposed as $\s=\hat{\s}\vee\check{\s}$ for
some $\hat{\s}\in\hat{\cC}(0)$ and $\check{\s}\in\check{\cC}(0)$. Note that $\hat{\cC}(0)$ has
$K=\ell_c(\ell_c-1)+2$ particles and $\check{\cC}(0)$ has $n_\b-K$ particles (and recall that, by
the first half of (\ref{LambdacondsKa}), $n_\b\to\infty$ as $\b\to\infty$). Define
\be{barc}
\cC^{\mathrm{fp}}(0) =  \big\{\s\in\tilde{\cC}(0)\colon\,
H_\b(\s)=H_\b(\hat{\s}) + H_\b(\check{\s})\big\},
\ee
i.e., the set of configurations consisting of a proto-critical droplet and a free particle
inside $B_{L_\b,L_\b}(0)$ not interacting with the particles outside $B_{L_\b,L_\b}(0)$.
Write $\cC^{\mathrm{fp},-}(0)$ and $\cC^{\mathrm{fp},+}(0)$ to denoting the subsets of
$\cC^{\mathrm{fp}}(0)$ where the free particle is at distance $L_\b$, respectively, 2 from
the proto-critical droplet. Split the right-hand side of (\ref{upkawII}) into a contribution
coming from  $\s,\s'\in\cC^{\mathrm{fp}}(0)$ and the rest:
\be{upkawIIsplit}
\textnormal{r.h.s.} (\ref{upkawII}) = |\L_\b|\,[II + \gamma_2(\b)],
\ee
where
\be{IIdef}
II = \min_{ {h\colon\,\cC^{\mathrm{fp}}(0) \to [0,1]}
\atop {h|_{\cC^{\mathrm{fp},-}(0)} \equiv 1,\,h|_{\cC^{\mathrm{fp},+}(0)}\equiv 0} }
\tfrac12 \sum_{\s,\s'\in\cC^{\mathrm{fp}}(0)}
\mu_\b(\s)c_\b(\s,\s')\,[h(\s)-h(\s')]^2
\ee
and $\gamma_2(\b)$ is an error term that will be estimated in Step 6. This term will turn
out to be small because of loss of entropy when the particle is at the boundary.

\medskip\noindent
{\bf 3.\ Reduction to capacity of simple random walk.} Estimate
\be{upkawIII}
\begin{aligned}
II
&= \min_{ {h\colon\,\cC^{\mathrm{fp}}(0) \to [0,1]}
\atop {h|_{\cC^{\mathrm{fp},-}(0)} \equiv 1,\,h|_{\cC^{\mathrm{fp},+}(0)} \equiv 0} } \\
&\qquad\qquad \tfrac12 \sum_{\check{\s},\check{\s}'\in\check{\cC}(0)}
\sum_{ {\hat{\s},\hat{\s}'\in\hat{\cC}(0):} \atop {\hat{\s}\vee\check{\s},\hat{\s}'
\vee\check{\s}'\in\cC^{\mathrm{fp}}(0)} }
\mu_\b(\hat{\s}\vee\check{\s})\,c_\b(\hat{\s}\vee\check{\s},\hat{\s}'\vee\check{\s}')
\,[h(\hat{\s}\vee\check{\s})-h(\hat{\s}'\vee\check{\s}')]^2\\
&\leq \min_{ {g\colon\,\hat{\cC}(0) \to [0,1]}
\atop {g|_{\hat{\cC}^-(0)} \equiv 1,\,g|_{\hat{\cC}^+(0)} \equiv 0} } \\
&\qquad\qquad \tfrac12 \sum_{\check{\s}\in\check{\cC}(0)}
\sum_{ {\hat{\s},\hat{\s}'\in\hat{\cC}(0):} \atop {\hat\s\vee\check{\s},\hat{\s}'
\vee\check{\s}\in\cC^{\mathrm{fp}}(0)} }
\mu_\b(\hat{\s}\vee\check{\s})\,c_\b(\hat{\s}\vee\check{\s},\hat{\s}'\vee\check{\s})
\,[g(\hat{\s})-g(\hat{\s}')]^2,
\end{aligned}
\ee
where $\hat{\cC}^-(0)$, $\hat{\cC}(0)^+$ denote the subsets of $\hat{\cC}(0)$ where
the free particle is at distance $L_\b$, respectively, 2 from the proto-critical
droplet, and the inequality comes from substituting
\be{hginsert}
h(\hat{\s}\vee\check{\s})= g(\hat{\s}), \qquad \hat{\s}\in\hat{\cC}(0),\,
\check{\s}\in\check{\cC}(0),
\ee
and afterwards replacing the double sum over $\check{\s},\check{\s}'\in\check{\cC}(0)$ by the
single sum over $\check{\s}\in\check{\cC}(0)$ because $c_\b(\hat{\s}\vee\check{\s},\hat{\s}'
\vee\check{\s}')>0$ only if either $\hat{\s}=\hat{\s}'$ or $\check{\s}=\check{\s}'$ (the
dynamics updates one site at a time). Next, estimate
\be{upkawIIIalt1}
\begin{aligned}
&\textnormal{r.h.s.} (\ref{upkawIII})\\
&\leq \sum_{\check{\s}\in\check{\cC}(0)} \frac{1}{Z_\b^{(n_\b)}}\,e^{-\b H_\b(\check{\s})}
\min_{ {g\colon\,\hat{\cC}(0) \to [0,1]}
\atop {g|_{\hat{\cC}^-(0)} \equiv 1,\,g|_{\hat{\cC}^+(0)}\equiv 0} }
\tfrac12 \sum_{ {\hat{\s},\hat{\s}'\in\hat{\cC}(0)} \atop
{\hat{\s}\vee\check{\s},\hat{\s}'\vee\check{\s}\in\cC^{\mathrm{fp}}(0)} } e^{-\b H_\b(\hat{\s})}
\,c_\b(\hat{\s},\hat{\s}')\,[g(\hat{\s})-g(\hat{\s}')]^2,
\end{aligned}
\ee
where we used $H_\b(\s)=H_\b(\hat{\s})+H_\b(\check{\s})$ from (\ref{barc}) and write
$c_\b(\hat{\s},\hat{\s}')$ to denote the transition rate associated with the Kawasaki
dynamics restricted to $B_{L_\b,L_\b}(0)$, which clearly equals $c_\b(\hat{\s}\vee
\check{\s},\hat{\s}'\vee\check{\s})$ for every $\check{\s}\in\check{\cC}(0)$ such
that $\hat{\s}\vee\check{\s},\hat{\s}'\vee\check{\s}\in\cC^{\mathrm{fp}}(0)$ because
there is no interaction between the particles inside and outside $B_{L_\b,L_\b}(0)$.
The minimum in the r.h.s.\ of (\ref{upkawIIIalt1}) can be estimated from above by
\be{Vbetasigmadef}
\sum_{\s\in\cP(0)} \cV_\b(\s)
\ee
with $\cP(0)$ the set of proto-critical droplets with lower-left corner at $0$, and
\be{Vdef}
\cV_\b(\s) = \min_{ {f\colon\,\Z^2 \to [0,1]}
\atop {f|_{P_\s(0)} \equiv 1,\,f|_{[B_{L_\b,L_\b}(0)]^c}\equiv 0} }
\tfrac12 \sum_{ {x,x'\in\Z^2}
\atop {x \sim x'} } [f(x)-f(x')]^2,
\ee
where $P_\s(0)$ is the support of the proto-critical droplet in $\s$, and $x \sim x'$ means
that $x$ and $x'$ are neighboring sites. Indeed, (\ref{Vbetasigmadef}) is obtained from the
expression in (\ref{upkawIIIalt1}) by dropping the restriction $\hat{\s}\vee\check{\s},
\hat{\s}'\vee\check{\s}\in\cC^{\mathrm{fp}}(0)$, substituting
\be{gfrel}
g(P_\s(0) \cup \{x\}) = f(x), \qquad \s\in\cP(0),\,x\in B_{L_\b,L_\b}(0)\backslash P_\s(0),
\ee
and noting that $c_\b(P_\s(0)\cup\{x\},P_\s(0)\cup\{x'\})=1$ when $x\sim x'$ and zero otherwise.
What (\ref{Vbetasigmadef}) says is that
\be{Vbetasigmadefalt}
\cV_\b(\s) = \CAPA(P_\s(0),[B_{L_\b.L_\b}(0)]^c)
\ee
is the capacity of simple random walk between the proto-critical droplet $P_\s(0)$ in $\s$ and
the exterior of $B_{L_\b.L_\b}(0)$. Now, define
\be{Cldef}
\check{Z}_\b^{(n-K)}(0) = \sum_{\check{\s}\in\check{\cC}(0)} e^{-\b H_\b(\check{\s})}.
\ee
Then we obtain via (\ref{Vbetasigmadef}) that
\be{upkawIIIalt2}
\textnormal{r.h.s.} (\ref{upkawIIIalt1})
\leq e^{-\b\Gamma^*}\,\frac{\check{Z}_\b^{(n-K)}(0)}{Z_\b^{(n_\b)}}
\sum_{\s\in\cP(0)} \cV_\b(\s),
\ee
where $\Gamma^*=-U[(\ell_c-1)^2+\ell_c(\ell_c-1)+1]$ is the binding energy of the proto-critical
droplet (compare with (\ref{Gammadef*})).

\medskip\noindent
{\bf 4.\ Capacity estimate.}
For future reference we state the following estimate on capacities for simple random walk.

\bl{capa.1}
Let $U\subset\Z^2$ be any set such that $\{0\} \subset U\subset B_{k,k}(0)$, with $k\in\N\cup\{0\}$
independent of $\b$. Let $V\subset\Z^2$ be any set such that $[B_{KL_\b,KL_\b}(0)]^c
\subset V \subset [B_{L_\b,L_\b}(0)]^c$, with $K\in\N$ independent of $\b$. Then
\be {capa.2}
\CAPA\left(\{0\},[B_{KL_\b,KL_\b}(0)]^c\right) \leq \CAPA\left(U,V\right)
\leq \CAPA\left(B_{k,k}(0),[B_{L_\b,L_\b}(0)]^c\right).
\ee
Moreover, via {\rm (\ref{Lbetadef}--\ref{dbprop})},
\be{capL}
\CAPA\left(B_{k,k}(0),[B_{KL_\b,KL_\b}(0)]^c\right)
= [1+o(1)]\,\frac{2\pi}{\log (KL_\b)-\log k}
= [1+o(1)]\,\frac{4\pi}{\b\Delta}, \qquad \b\to\infty.
\ee
\el

\bpr
The inequalities in \eqref{capa.2} follow from standard monotonicity properties of
capacities. The asymptotic estimate in (\ref{capL}) for capacities of concentric
boxes are standard (see e.g.\ Lawler~\cite{L91}, Section 2.3), and also follow by
comparison to Brownian motion.
\epr

We can apply Lemma~\ref{capa.1} to estimate $\cV_\b(\sigma)$ in (\ref{Vbetasigmadefalt}), since
the proto-critical droplet with lower-left corner in 0 fits inside the box $B_{2\ell_c,2\ell_c}(0)$.
This gives
\be{Vrelcap}
\cV_\b(\s) = \frac{4\pi}{\b\Delta}\,[1+o(1)], \qquad \forall\,\s\in\cP(0),\,\b\to\infty.
\ee
Morover, from Bovier, den Hollander, and Nardi~\cite{BdHN06}, Lemmas 3.4.2--3.4.3, we know
that $N=|\cP(0)|$, the number of shapes of the proto-critical droplet, equals $N=\frac13
\ell_c^2(\ell_c^2-1)$.

\medskip\noindent
{\bf 5.\ Equivalence of ensembles.}
According to Lemma~\ref{cgcequiv} in Appendix~\ref{appB}, we have
\be{muCnrellim}
\frac{\check{Z}_\b^{(n_\b-K)}(0)}{Z_\b^{(n_\b)}}
= (\rho_\b)^K\,\mu_\b(\cS)\,[1+o(1)], \qquad \b\to\infty.
\ee
This is an ``equivalence of ensembles'' property relating the probabilities to
find $n_\b-K$, respectively, $n_\b$ particles inside $[B_{L_\b,L_\b}(0)]^c$ (recall
(\ref{cchec})). Combining (\ref{upkaw}--\ref{upkawsplit}), (\ref{upkawII}),
(\ref{upkawIIsplit}), (\ref{upkawIII}), (\ref{upkawIIIalt1}), (\ref{upkawIIIalt2})
and (\ref{Vrelcap}--\ref{muCnrellim}), we get
\be{Iestfin}
\CAPA(\cS,\cC^+) \leq \gamma_1(\b) + |\L_\b|\gamma_2(\b)
+ N\,|\L_\b|\,\frac{4\pi}{\b\D}\,e^{-\b\G}\,\mu_\b(\cS)\,[1+o(1)],
\qquad \b\to\infty,
\ee
where we use that $\G^*+\D K=\G$ defined in (\ref{Gammadef*}). This completes the proof
of the upper bound, provided that the error terms $\gamma_1(\b)$ and $\gamma_2(\b)$ are
negligible.

\medskip\noindent
{\bf 6.\ Second error term.} To estimate the error term $\gamma_2(\b)$, note that the configurations
in $\tilde{\cC}(0)\backslash\cC^{\mathrm{fp}}(0)$ are those for which inside $B_{L_\b,L_\b}
(0)$ there is a proto-critical droplet whose lower-left corner is at 0, and a particle
that is at the boundary and attached to some cluster outside $B_{L_\b,L_\b}(0)$.
Recalling (\ref{upkawII}--\ref{IIdef}), we therefore have
\be{gamma2est}
\gamma_2(\b) \leq
\sum_{\s\in \tilde{\cC}(0)\backslash\cC^\mathrm{fp}(0)}\,\,\sum_{\s'\in \tilde{\cC}(0)}
\mu_\b(\s)c_\b(\s,\s')\,[h(\s)-h(\s')]^2
\leq 6 \mu_\b(\tilde{\cC}(0)\backslash\cC^{\mathrm{fp}}(0)),
\ee
where we use that $h\colon\tilde{\cC}(0)\to[0,1]$, $\mu_\b(\s)c_\b(\s,\s') = \mu_\b(\s)
\wedge \mu_\b(\s')$, and there are 6 possible transitions from $\tilde{\cC}(0)\backslash
\cC^{\mathrm{fp}}(0)$ to $\tilde{\cC}(0)$: 3 through a move by the particle at the
boundary of $B_{L_\b,L_\b}(0)$ and 3 through a move by a particle in the cluster outside
$B_{L_\b,L_\b}(0)$. Since
\be{ener1}
H_\b(\s) \geq H_\b(\hat{\s}) + H_\b(\check{\s}) - U,
\qquad \s\in\tilde{\cC}(0)\backslash\cC^{\mathrm{fp}}(0),
\ee
it follows from the same argument as in Steps 3 and 5 that
\be{gamma2}
\mu_\b(\tilde{\cC}(0)\backslash\cC^{\mathrm{fp}}(0))
\leq N\,e^{-\b\Gamma^*}\,(\rho_\b)^{K+1}\,\mu_\b(\cS)\,
e^{\b U}\,4(K-1)\,[1+o(1)],
\ee
where $(\rho_\b)^{K+1}$ comes from the fact that $n_\b-(K+1)$ particles are outside
$B_{L_\b-1,L_\b-1}(0)$ (once more use Lemma~\ref{cgcequiv} in Appendix~\ref{appB}),
$e^{\b U}$ comes from the gap in (\ref{ener1}), and $4(K-1)$ counts the maximal
number of places at the boundary of $B_{L_\b,L_\b}(0)$ where the particle can
interact with particles outside $B_{L_\b,L_\b}(0)$ due to the constraint that
defines $\cS$ (recall Definition~\ref{Rdef*})(a)). Since $\rho_\b e^{\b U}=o(1)$
by (\ref{metreg*}), we therefore see that $\gamma_2(\b)$ indeed is small compared
to the main term of (\ref{Iestfin}).

\medskip\noindent
{\bf 7.\ First error term.} To estimate the error term $\gamma_1(\b)$, we define the sets of
pairs of configurations
\be{I1I2}
\begin{aligned}
\cI_1 &= \{(\s,\eta)\in [\cX_\b^{(n_\b)}]^2\colon\,
\s\in\cS,\eta\in \cS^c\backslash\tilde{\cC}\},\\
\cI_2 &= \{(\s,\eta)\in [\cX_\b^{(n_\b)}]^2\colon\,
\s\in\tilde{\cC},\eta\in \cS^c\backslash\tilde{\cC}\},
\end{aligned}
\ee
and estimate
\be{gamma12u}
\gamma_1(\b)\leq \tfrac12 \sum_{i=1}^2\sum_{(\s,\eta)\in \cI_i}\mu_\b(\s)\,c_\b(\s,\eta)
 = \tfrac12\Sigma(\cI_1) + \tfrac12\Sigma(\cI_2).
\ee
The sum $\Sigma(\cI_1)$ can be written as
\be{I1u}
\Sigma(\cI_1) = |\L_\b| \sum_{\s\in \cP} \sum_{\eta\in\cS^c\backslash\tilde{\cC}}
c_\b(\eta,\s)\,1\Big\{|\textnormal{supp}[\eta]\cap B_{L_\b,L_\b}(0)|=K\Big\}
\,\frac{1}{Z_\b^{(n_\b)}}\,e^{-\b H_\beta(\eta)},
\ee
where we use that $\mu_\b(\s)c_\b(\s,\eta)=\mu_\b(\eta)c_\b(\eta,\s)$, $\s,\eta\in
\cX_\b^{(n_\b)}$, and $c_\b(\eta,\s)=0$, $\eta\in\cS^c\backslash\tilde{\cC}$, $\s\notin\cP$
(recall Definition~\ref{Rdef*}(b)). We have
\be{Hest12}
H_\beta(\eta) \geq H_\b(\hat{\eta}) + H_\b(\check{\eta})
-kU, \qquad \eta \in \cS^c\backslash\tilde{\cC},
\ee
where $k$ counts the number of pairs of particles interacting across the boundary of
$B_{L_\b,L_\b}(0)$. Moreover, since $\eta\notin \tilde{C}$, we have
\be{Hest13}
H_\b(\hat{\eta}) \geq \Gamma^* + U.
\ee
Inserting (\ref{Hest12}--\ref{Hest13}) into (\ref{I1u}), we obtain
\be{I1u2}
\begin{aligned}
\Sigma(\cI_1)
&\leq |\L_\b|\,e^{-\b\Gamma^*}\,
\mu_\b(\cS)\,[1+o(1)] \sum_{k=0}^K (\rho_\b)^{K+k}\,[4(K-1)]^k\,e^{\b(k-1)U}\\
&= |\L_\b|\,e^{-\b\Gamma}
\mu_\b(\cS)\,[1+o(1)]\,e^{-\b U},
\end{aligned}
\ee
where $(\rho_\b)^{K+k}$ comes from the fact that $n_\b-(K+k)$ particles are outside
$B_{L_\b-1,L_\b-1}(0)$ (once more use Lemma~\ref{cgcequiv} in Appendix~\ref{appB}), and
the inequality again uses an argument similar as in Steps 3 and 5. Therefore
$\Sigma(\cI_1)$ is small compared to the main term of (\ref{Iestfin}). The sum
$\Sigma(\cI_2)$ can be estimated as
\be{i2}
\begin{aligned}
\Sigma(\cI_2)
&= \sum_{\s\in\tilde{\cC}} \sum_{\eta\in \cS^c\backslash\tilde{\cC}}
\mu_\b(\s)\,c_\b(\s,\eta)\\
&= |\L_\b| \sum_{\s\in\tilde{\cC}(0)} \mu_\b(\s)
\sum_{\eta\in \cS^c\backslash\tilde{\cC}(0)} c_\b(\s,\eta)\\
&\leq |\L_\b|\,\mu_\b(\tilde{\cC}(0))\,\big\{e^{-\b\,U} + (4L_\b)\,\rho_\b\,[1+o(1)]\big\},
\end{aligned}
\ee
where the first term comes from detaching a particle from the critical droplet and the
second term from a extra particle entering $B_{L_\b,L_\b}(0)$. The term between braces
is $o(1)$. Moreover, $\mu_\b(\tilde{\cC}(0))=\mu_\b(\cC^\mathrm{fp}(0))+ \mu_\b(\tilde{\cC}
(0)\backslash\cC^\mathrm{fp}(0))$. The second term was estimated in (\ref{gamma2}), the
first term can again be estimated as in Steps 3 and 5:
\be{muCcal}
\mu_\b(\cC^\mathrm{fp}(0))
= \sum_{\hat{\s}\in\hat{\cC}(0)}
\sum_{ {\check{\s}\in\check{\cC}(0)} \atop {\hat{\s}\vee\check{\s}\in\cC^\mathrm{fp}(0)} }
\mu_\b(\hat{\s}\vee\check{\s})
= N\,e^{-\b\G^*}\,\frac{\check{Z}_\b^{(n_\b-K)}(0)}{Z_\b^{(n_\b)}}
= N\,e^{-\b\G}\,\mu_\b(\cS)\,[1+o(1)].
\ee
Therefore also $\Sigma(\cI_2)$ is small compared to the main term of (\ref{Iestfin}).

\bigskip\noindent
\underline{Lower bound}:
The proof of the lower bound follows the same line of argument as for Glauber dynamics
in that it relies on the construction of a suitable unit flow. This flow will, however,
be considerably more difficult. In particular, we will no longer be able to get away
with choosing a deterministic flow, and the full power of the Berman-Konsowa variational
principle has to be brought to bear. The proof comes in 5 steps.

For future reference we state the following property of the harmonic function for simple
random walk on $\Z^2$.

\bl{flow-lb.11}
Let $g$ be the harmonic function of simple random walk on $B_{2L_\b,2L_\b}(0)$ (which
is equal to $1$ on $\{0\}$ and $0$ on $[B_{2L_\b,2L_\b}(0)]^c$). Then there exists a
constant $C<\infty$ such that
\be{flow-nb.12}
\sum_{e}[g(z)-g(z+e)]_+ \leq C/L_\b \qquad \forall\,z\in [B_{L_\b,L_\b}(0)]^c.
\ee
\el

\bpr
See e.g.\ Lawler, Schramm and Werner~\cite{LSW04}, Lemma 5.1. The proof can be given
via the estimates in Lawler~\cite{L91}, Section 1.7, or via a coupling argument.
\epr

\medskip\noindent
{\bf 1.\ Starting configurations.}
We start our flow on a subset of the configurations in $\cC^+$ that is
sufficiently large and sufficiently convenient. Let $\cC_2^+\subset\cC^+$
denote the set of configurations having a proto-critical with lower-left corner
at some site $x\in\L_\b$, a free particle at distance 2 from this proto-critical
droplet, no other particles in the box $B_{2L_\b,2L_\b}(x)$, and satisfying the
constraints in $\cS_L$, i.e., all other boxes of size $2L_\b$ carry no more particles
than there are in a proto-critical droplet. This is the same as $\cC^+$, except that
the box around the proto-critical droplet has size $2L_\b$ rather than $L_\b$.

Let $K=\ell_c(\ell_c-1)+2$ be the volume of the critical droplet, and let $\cS_2^{(n_\b-K)}$
be the analogue of $\cS$ when the total number of particles is $n_\b-K$ and the boxes in
which we count particles have size $2L_\b$ (compare with Definition~\ref{Rdef*}). Similarly
as in (\ref{lb2}), our task is to derive a lower bound for $\CAPA(\cS_L,(\cS^c\backslash
\tilde{\cC})\cup\cC^+) = \CAPA((\cS^c\backslash\tilde{\cC})\cup\cC^+,\cS_L) \geq
\CAPA(\cC_L,\cS_L)$, where $\cC_L\subset\cC_2^+\subset\cC^+$ defined by
\be{DLdefKaw}
\cC_L = \{\s \cup P_{(y)}(x,z)\colon\,\s\in\cS_2^{(n_\b-K)},\,
x,y\in\L_\b,\,(x,y,z) \bot \s\}
\ee
is the analog of (\ref{rcap}), namely, the set of configurations obtained from
$\cS_2^{(n_\b-K)}$ by adding a critical droplet somewhere in $\L_\b$ (lower-left
corner at $x$, protuberance at $y$, free particle at $z$) such that it does not
interact with the particles in $\s$ and has an empty box of size $2L_\b$ around it.
Note that the $n_\b-K$ particles can block at most $n_\b(2L_\b)^2=o(|\L_\b|)$ sites
from being the center of an empty box of size $2L_\b$, and so the critical particle
can be added at $|\L_\b|-o(|\L_\b|)$ locations.

We partition $\cC_L$ into sets $\cC_L(x)$, $x\in\L_\b$, according to the location of
the proto-critical droplet. It suffices to consider the case where the critical droplet
is added at $x=0$, because the union over $x$ trivially produces a factor $|\L_\b|$.

\medskip\noindent
{\bf 2.\ Overall strategy.}
Starting from a configuration in $\cC_L(0)$, we will successively pick $K-L$ particles
from the critical droplet (starting with the free particle at $z$ at distance $2$) and
move them out of the box $B_{L_\b,L_\b}(0)$, placing them essentially uniformly in the
annulus $B_{2L_\b,2L_\b}(0) \backslash B_{L_\b,L_\b}(0)$. Once this has been achieved,
the configuration is in $\cS_L$. Each such move will produce an entropy of order $L_\b^2$,
which will be enough to compensate for the loss of energy in tearing down the droplet
(recall Fig.~\ref{fig:nuclpath*}). The order in which the particles are removed
follows the canonical order employed in the lower bound for Glauber dynamics
(recall Fig.~\ref{fig:droptear}). As for Glauber, we will use Proposition~\ref{capMarkov}
to estimate
\be{BKflowKawa}
\CAPA(\cC_L,\cS_L) \geq |\L_\b| \sum_{\s\in\cC_L(0)}
\sum_{\g\colon\,\g_0=\s} \P^f(\g) \sum_{k=0}^{\tau(\g)}
\left[\frac{f(\g_k,\g_{k+1})}{\mu_\b(\g_k)c_\b(\g_k,\g_{k+1})}\right]^{-1}
\ee
for a suitably constructed flow $f$ and associated path measure $\P^f$, starting from
some initial distribution on $\cC_L(0)$ (which as for Glauber will be irrelevant), and
$\tau(\g)$ the time at which the last of the $K-L$ particles exits the box
$B_{L_\b,L_\b}(0)$ .

The difference between Glauber and Kawasaki is that, while in Glauber the droplet
can be torn down via single spin-flips, in Kawasaki after we have detached a particle
from the droplet we need to move it out of the box $B_{L_\b,L_\b}(0)$, which takes
a large number of steps. Thus, $\tau(\gamma)$ is the sum of $K-L$ stopping times,
each except the first of which is a sum of two stopping times itself, one to detach
the particle and one to move it out of the box $B_{L_\b,L_\b}(0)$. With each motion
of a single particle we need to gain an entropy factor of order close to $1/\rho_\b$.
This will be done by constructing a flow that involves only the motion of this single
particle, based on the harmonic function of the simple random walk in the box
$B_{2L_\b,2L_\b}(0)$ up to the boundary of the box $B_{L_\b,L_\b}(0)$. Outside
$B_{L_\b,L_\b}(0)$ the flow becomes more complex: we modify it in such a way that
a small fraction of the flow, of order $L_\b^{-1+\e}$ for some $\e>0$ small enough, is
going into the direction of removing the next particle from the droplet. The reason for this
choice is that we want to make sure that the flow becomes sufficiently small, of order
$L_\b^{-2+\e}$, so that this can compensate for the fact that the Gibbs weight in the
denominator of the lower bound in (\ref{ab.5}) is reduced by a factor $e^{-\b U}$ when
the protuberance is detached. The reason for the extra $\e$ is that we want to make sure
that, along most of the paths, the protuberance is detached before the first particle
leaves the box $B_{2L_\b,2L_\b}(0)$.

Once the protuberance detaches itself from the proto-critical, the first particle
stops and the second particle moves in the same way as the first particle did when
it moved away from the proto-critical droplet, and so on. This is repeated until no
more than $L$ particles remain in $B_{L_\b,L_\b}(0)$, by which time we have reached
$\cS_L$. As we will see, the only significant contribution to the lower bound comes
from the motion of the first particle (as for Glauber), and this coincides with the
upper bound established earlier. The details of the construction are to some extent
arbitrary and there are many other choices imaginable.

\medskip\noindent
{\bf 3.\ First particle.}
We first construct the flow that moves the particle at distance $2$ from the
proto-critical droplet to the boundary of the box $B_{L_\b,L_\b}(0)$. This flow
will consist of independent flows for each fixed shape and location of the critical
droplet. This first part of the flow will be seen to produce the essential contribution
to the lower bound.

We label the configurations in $\cC_L(0)$ by $\s$, describing the shape of the critical
droplet, as well as the configuration outside the box $B_{2L_\b,2L_\b}(0)$, and we label
the position of the free particle in $\s$ by $z_1(\s)$.

Let $g$ be the harmonic function for simple random walk with boundary
conditions $0$ on $[B_{2L_\b,2L_\b}(0)]^c$ and $1$ on the critical droplet.
Then we choose our flow to be
\be{flow-lb.1}
f(\s(z),\s(z'))
= \begin{cases}
C_1\,[g(z)-g(z+e)]_+, &\mbox{if } z'=z+e, \|e\|=1,\\
0,                    &\mbox{otherwise},
\end{cases}
\ee
where $\s(z)$ is the configuration obtained from $\s$ by placing the first particle at
site $z$. The constant $C_1$ is chosen to ensure that $f$ defines a unit flow in the
sense of Definition~\ref{unit-flow}, i.e.,
\be{flow-lb.2}
\sum_{\s\in \cC_L(0)} C_1 \sum_{z_1(\s),e}
[g(z_1(\s))-g(z_1(\s)+e)]
= C_1\,\sum_{\s\in\cC_L(0)}\,\CAPA \left(P_\s(0),[B_{2L_\b,2L_\b}(0)]^c\right)
= 1,
\ee
where $P_\s(0)$ denotes the support of the proto-critical droplet in $\sigma$, and the
capacity refers to the simple random walk.

Now, let $z^1(k)$ be the location of the first particle at time $k$, and
\be{taufirstdef}
\t^1 = \inf\{k\in\N\colon\,z^1(k) \in [B_{L_\b,L_\b}(0)]^c\}
\ee
be the first time when, under the Markov chain associated to the flow $f$, it
exits $B_{L_\b,L_\b}(0)$. Let $\g$ be a path of this Markov chain. Then, by
(\ref{flow-lb.1}--\ref{flow-lb.2}), we have
\be{flow-lb.3}
\sum_{k=0}^{\t^1} \frac {f(\g_k,\g_{k+1})}{\mu_\b(\g_k) c_\b(\g_k,\g_{k+1})}
= \frac{C_1[g(z^1(0))-g(z^1(\t^1))]}{\mu_\b(\g_0)}
\ee
where the sum over the $g$'s is telescoping because only paths along which the
$g$-function decreases carry positive probability, and $c_\b(\g_k,\g_{k+1})=1$
for all $0\leq k \leq \t^1$ because the first particle is free. We have
$g(z^1(0))=1$, while, by Lemma~\ref{flow-lb.11}, there exists a $C<\infty$ such
that
\be{gxest}
g(x) \leq C/\log L_\b, \qquad x \in [B_{L_\b,L_\b}(0)]^c.
\ee
Therefore
\be{flow-lb.3rew}
\sum_{k=0}^{\t^1} \frac {f(\g_k,\g_{k+1})}{\mu_\b(\g_k) c_\b(\g_k,\g_{k+1})}
= \frac{C_1}{\mu_\b(\g_0)}\,[1+o(1)].
\ee
Next, by Lemma~\ref{capa.1}, we have
\be{CapV}
\CAPA\big(P_\s(0),[B_{2L_\b,2L_\b}(0)]^c\big) = \frac{4\pi}{\b\D}\,[1+o(1)],
\qquad \s\in\cC_L(0),\,\b\to\infty,
\ee
(because $\{0\} \subset P_\s(0) \subset B_{2\ell_c.2\ell_c}(0)$ for all $\s\in\cC_L(0)$).
Since $N=|\cC_L(0)|$, it follows from (\ref{flow-lb.2}) that
\be{C1asymp}
\frac{1}{C_1} = N\,\frac{4\pi}{\b\D}\,[1+o(1)],
\ee
and so (\ref{flow-lb.3rew}) becomes
\be{flow-lb.4}
\left[\sum_{k=0}^{\t^1} \frac {f(\g_k,\g_{k+1})}{\mu_\b(\g_k) c_\b(\g_k,\g_{k+1})}\right]^{-1}
= \mu_\b(\g_0)\,N\,\frac{4\pi}{\b\D}\,[1+o(1)].
\qquad \b\to\infty,
\ee
This is the contribution we want, because when we sum (\ref{flow-lb.4}) over $\g_0=\s\in\cC_L(0)$
(recall (\ref{BKflowKawa})), we get a factor
\be{Kawlbrels}
\mu_\b(\cC_L(0)) = e^{-\b\G}\,\mu_\b(\cS)\,[1+o(1)].
\ee
To see why (\ref{Kawlbrels}) is true, recall from (\ref{DLdefKaw}) that $\cC_L(0)$ is obtained
from $\cS_2^{(n_\b-K)}$ by adding a critical droplet with lower-left corner at the origin
that does not interact with the $n_\b-K$ particles elsewhere in $\L_\b$. Hence
\be{Kawlbrels2}
\mu_\b(\cC_L(0)) = e^{-\b\Gamma^*}\,
\frac{\tilde{Z}_\b^{(n_\b-K)}(0)}{Z_\b^{(n_\b)}},
\ee
where $\tilde{Z}_\b^{(n_\b-K)}(0)$ is the analog of $\check{Z}_\b^{(n_\b-K)}(0)$ (defined in
(\ref{Cldef})) obtained by requiring that the $n_\b-K$ particles are in $[R_{\ell_c,\ell_c}(0)]^c$
instead of $[B_{L_\b,L_\b}(0)]^c$. However, it will follow from the proofs of
Lemmas~\ref{cgcequiv}--\ref{checkSrel} in Appendix~\ref{appB} that, just as in (\ref{muCnrellim}),
\be{Kawlbrels2alt}
\frac{\tilde{Z}_\b^{(n_\b-K)}(0)}{Z_\b^{(n_\b)}} = (\rho_\b)^K\,\mu_\b(\cS)\,[1+o(1)],
\qquad \b\t\infty,
\ee
which yields (\ref{Kawlbrels}) because $\Gamma=\Gamma^*+K\Delta$. For the remaining part of the
construction of the flow it therefore suffices to ensure that the sum beyond $\t^1$ gives
a smaller contribution.

\medskip\noindent
{\bf 4.\ Second particle.}
Once the first particle (i.e., the free particle) has left the box $B_{L_\b,L_\b}(0)$,
we need to allow the second particle (i.e., the protuberance) to detach itself from
the proto-critical droplet and to move out of $B_{L_\b,L_\b}(0)$ as well. The
problem is that detaching the second particle reduces the Gibbs weight appearing
in the denominator by $e^{-U\b}$, while the increments of the flow are reduced
only to about $1/L_\b$. Thus, we cannot immediately detach the second particle.
Instead, we do this with probability $L_\b^{-1+\e}$ only.

The idea is that, once the first particle is outside $B_{L_\b,L_\b}(0)$, we
\emph{leak} some of the flow that drives the motion of the first particle into
a flow that detaches the second particle. To do this, we have to first construct
a \emph{leaky flow} in $B_{2L_\b,2L_\b}(0) \backslash B_{L_\b,L_\b}(0)$ for
simple random walk. This goes as follows.

Let $p(z,z+e)$ denote the transition probabilities of simple random walk driven
by the harmonic function $g$ on $B_{2L_\b,2L_\b}(0)$. Put
\be{tilpdef}
\tilde{p}(z,z+e) = \left\{\begin{array}{ll}
p(z,z+e), &\mbox{if } z \in B_{L_\b,L_\b}(0),\\
(1-L_\b^{-1+\e})\, p(z,z+e),
&\mbox{if } z \in B_{2L_\b,2L_\b}(0) \backslash B_{L_\b,L_\b}(0).
\end{array}
\right.
\ee
Use the transition probabilities $\tilde{p}(z,z+e)$ to define a path measure $\tilde{P}$.
This path measure describes simple random walk driven by $g$, but with a killing probability
$L_\b^{-1+\e}$ inside the annulus $B_{2L_\b,2L_\b}(0) \backslash B_{L_\b,L_\b}(0)$. Put
\be{kleakdef}
k(z,z+e) = \sum_{\g} \tilde{P}(\g) \1_{(z,z+e) \in \g},
\qquad z \in B_{2L_\b,2L_\b}(0).
\ee
This edge function satisfies the following equations:
\be{flow-lb.6}
\begin{aligned}
&\bullet k(z,z+e) = [g(z)-g(z+e)]_+,\\
&\qquad \mbox{if } z \in B_{L_\b,L_\b}(0),\\
&\bullet k(z,z+e) = 0,\\
&\qquad \mbox{if } z \in B_{2L_\b,2L_\b}(0)\backslash B_{L_\b,L_\b}(0)
\mbox{ and } [g(z)-g(z+e)]_+=0,\\
&\bullet (1-L_\b^{-1+\e}) \sum_{e} k(z+e,z) \1_{g(z+e)-g(z)>0}
= \sum_{e} k(z,z+e)\1_{g(z) -g(z+e)>0}\\
&\qquad \mbox{if }\ z\in B_{2L_\b,2L_\b}(0)\backslash B_{L_\b,L_\b}(0).
\end{aligned}
\ee
Note that inside the annulus $B_{2L_\b,2L_\b}(0)\backslash B_{L_\b,L_\b}(0)$ at each site
the flow out is less than the flow in by a leaking factor $1-L_\b^{-1+\e}$. We pick $\e>0$
so small that
\be{epscond}
e^{\b U} \mbox{ is exponentially smaller in } \b \mbox{ than } L_\b^{2-\e},
\ee
(which is possible by (\ref{metreg*}) and (\ref{Lbetadef}--\ref{dbprop})). The important
fact for us is that this leaky flow is dominated by the harmonic flow associated with $g$,
in particular, the flow in satisfies
\be{flow-lb.7}
\sum_{e} k(z+e,z) \leq \sum_{e} [g(z+e)-g(z)]_+ \qquad
\forall\,z\in  B_{2L_\b,2L_\b}(0),
\ee
(and the same applies for the flow out). This inequality holds because $g$ satisfies the
same equations as in (\ref{tilpdef}--\ref{kleakdef}) but without the leaking factor
$1-L_\b^{-1+\e}$.

Using this leaky flow, we can now construct a flow involving the first two
particles, as follows:
\be{flow-lb.8}
\begin{aligned}
&\bullet f(\s(z_1,a),\s(z_1+e,a)) = C_1 k(z_1,z_1+e),\\
&\qquad \mbox{if } z_1\in B_{2L_\b,2L_\b}(0),\\
&\bullet f(\s(z_1,a),\s(z_1,b)) = C_1 L_\b^{-1+\e} \sum_{e} k(z_1,z_1+e),\\
&\qquad \mbox{if } z_1\in B_{2L_\b,2L_\b}(0) \backslash B_{L_\b,L_\b}(0),\\
&\bullet f(\s(z_1,z_2),\s(z_1,z_2+e)) = \left\{C_1 L_\b^{-1+\e}
\sum_{e} k(z_1,z_1+e)\right\} [g(z_2)- g(z_2+e)]_+,\\
&\qquad \mbox{if } z_1\in B_{2L_\b,2L_\b}(0) \backslash B_{L_\b,L_\b}(0),
z_2 \in B_{L_\b,L_\b}(0)\backslash P_\s(0).
\end{aligned}
\ee
Here, we write $a$ and $b$ for the locations of the second particle prior and
after it detaches itself from the proto-critical droplet, and $\s(z_1,z_2)$ for
the configuration obtained from $\s$ by placing the first particle (that was at
distance $2$ from the proto-critical droplet) at site $z_1$ and the second particle
(that was the protuberance) at site $z_2$. The flow for other motions is zero, and
the constant $C_1$ is the \emph{same} as in (\ref{flow-lb.1}--\ref{flow-lb.2})

We next define two further stopping times, namely,
\be{flow-lb.9}
\zeta^2 = \inf\{k\in\N\colon\,z^2(\g_k)=b\},
\ee
i.e., the first time the second particle (the protuberance) detaches itself
from the proto-critical droplet, and
\be{flow-lb.10}
\t^2 = \inf\{k\in\N\colon\,z^2(\g_k)\in [B_{L_\b,L_\b}(0)]^c\},
\ee
i.e., the first time the second particle exits the box $B_{L_\b,L_\b}(0)$.
Note that, since we choose the leaking probability to be $L^{-1+\e}$, the
probability that $\zeta^2$ is larger than the first time the first particle
exits $B_{2L_\b,2L_\b}(0)$ is of order $\exp[-L_\b^{\epsilon}]$ and hence
is negligible. We will disregard the contributions of such paths in the lower
bound. Paths with this property will be called \emph{good}.

We will next show that (\ref{flow-lb.3}) also holds if we extend the sum along
any path of positive probability up to $\zeta^2$. The reason for this lies in
Lemma~{flow-lb.11}. Let $\g$ be a path that has a positive probability under
the path measure $\P^{f}$ associated with $f$ stopped at $\t^2$. We will assume
that this path is good in the sense described above. To that end we decompose
\be{flow-lb.13}
\begin{aligned}
&\sum_{k=0}^{\t^2} \frac {f(\g_k,\g_{k+1})}{\mu_\b(\g_k) c_\b(\g_k,\g_{k+1})}\\
&= \sum_{k=0}^{\t^1} \frac {f(\g_k,\g_{k+1})}{\mu_\b(\g_k) c_\b(\g_k,\g_{k+1})}
+ \sum_{k=\t^1+1}^{\zeta^2-2} \frac {f(\g_k,\g_{k+1})}{\mu_\b(\g_k) c_\b(\g_k,\g_{k+1})}
+\sum_{k=\zeta^2-1}^{\t^2} \frac {f(\g_k,\g_{k+1})}{\mu_\b(\g_k)c_\b(\g_k,\g_{k+1})}\\
&= I+II+III.
\end{aligned}
\ee
The term $I$ was already estimated in (\ref{flow-lb.3}--\ref{Kawlbrels}). To estimate $II$,
we use (\ref{gxest}) and (\ref{flow-lb.7}--\ref{flow-lb.8}) to bound (compare with
(\ref{flow-lb.3}))
\be{flow-lb.13.1}
II \leq
C_1\,\frac{g(z^1(\zeta^2))-g(z^1(\t^1))}{\mu_\b(\g_0)}
\leq  C_1\,\frac {[C/\log L_\b]}{\mu_\b(\g_0)},
\ee
which is negligible compared to $I$ due to the factor $C/\log L_\b$. It remains to
estimate $III$. Note that
\be{flow-lb.14}
III = \frac{f(\g_{\zeta^2-1},\g_{\zeta^2})}{\mu_\b(\g_{\zeta^2-1})
c_\b(\g_{\zeta^2-1},\g_{\zeta^2})}
+\sum_{k=\zeta^2}^{\t^2} \frac{f(\g_k,\g_{k+1})}{\mu_\b(\g_k)c_\b(\g_k,\g_{k+1})}.
\ee
The first term corresponds to the move when the protuberance detaches itself from the
proto-critical droplet. Its numerator is given by $f(\s(z_1,a),\s(z_1,b))$ (for some
$z_1\in [B_{L_\b,L_\b}(0)]^c$) which, by Lemma \ref{flow-lb.11} and
(\ref{flow-lb.7}--\ref{flow-lb.8}), is smaller than $C_1L_\b^{-1+\e}CL_\b^{-1}
= C_1CL_\b^{-2+\e}$. On the other hand, its denominator is given by
\be{flow-lb.14.13}
\mu(\g_{\zeta^2-1}) c_\b(\g_{\zeta^2-1},\g_{\zeta^2}) =\mu_\b(\g_0) e^{-U\b}.
\ee
The same holds for the denominators in all the other terms in $III$, while the
numerators in these terms satisfy the bound
\be{flow-nb.14.4}
f(\g_k,\g_{k+1}) \leq C_1\,C\,L_\b^{-2+\e} \big[g(z^2(\g_k))-g(z^2(\g_{k+1}))\big].
\ee
Adding up the various terms, we get that
\be{flow-lb.14.5}
III \leq \frac{C_1}{\mu_\b(\g_0)} L_\b^{-2+\e}e^{\b U}
\Big(1 + [g(z^2(\zeta^2))-g(z^2(\t^2)]\Big)
\leq \frac{2C_1}{\mu_\b(\g_0)} L_\b^{-2+\e}e^{\b U}.
\ee
The right-hand side is smaller than $I$ by a factor $L_\b^{-2+\e} e^{\b U}$,
which, by (\ref{epscond}), is exponentially small in $\b$.

\medskip\noindent
{\bf 5.\ Remaining  particles.}
The lesson from the previous steps is that we can construct a flow with
the property that each time we remove a particle from the droplet we gain
a factor $L_\b^{-2+\e}$, i.e., almost $e^{-\D\b}$. (This entropy gain corresponds
to the gain from the magnetic field in Glauber dynamics, or from the activity
in Kawasaki dynamics on a finite open box.) We can continue our flow by tearing
down the critical droplet in the same order as we did for Glauber dynamics. Each
removal corresponds to a flow that is built in the same way as described in Step 4
for the second particle. There will be some minor modifications involving a negligible
fraction of paths where a particle hits a particle that was moved out earlier, but
this is of no consequence. As a result of the construction, the sums along the
remainders of these paths will give only negligible contributions.

Thus, we have shown that the lower bound coincides, up to a factor $1+o(1)$,
with the upper bound and the lemma is proven.
\epr

\subsection{Proof of Theorem \ref{HNKawasaki}(b)}
\label{S4.2}

The same observation holds as in (\ref{eqpotMdrop}).

\bpr
The proof of Theorem~\ref{HNKawasaki}(b) follows along the same lines as that of
Theorem~\ref{HNKawasaki}(a). The main point is to prove that $\CAPA(\cD_M,\cS_L)
=[1+o(1)]\CAPA(\cC^+,\cS_L)$. Since $\CAPA(\cS_L,\cD_M)\leq\CAPA(\cS_L,\cC^+)$, we
need only prove a lower bound on $\CAPA(\cD_M,\cS_L)$. This is done in almost
exactly the same way as for Glauber, by using the construction given there and
substituting each Glauber move by a flow involving the motion of just two particles.

Note that, as long as $M=e^{o(\beta)}$, an $M \times M$ droplet can be added
at $|\L_\b|-o(|\L_\b|)$ locations to a configuration $\s\in\cS$ (compare with
(\ref{DLdefKaw})). The only novelty is that we have to eventually remove the
cloud of particles that is produced in the annulus $B_{2L_\b,2L_\b}(0) \backslash
B_{L_\b,L_\b}(0)$. This is done in much the same way as before. As long as only
$e^{o(\b)}$ particles have to be removed, potential collisions between particles
can be ignored as they are sufficiently unlikely.
\epr


\appendix


\section{Appendix: sparseness of subcritical droplets}
\label{appA}
\setcounter{equation}{0}

Recall Definition~\ref{Rdef}(a) and (\ref{fprop}--\ref{Wdef}). In this section we
prove (\ref{Wasymp}).

\bl{wses}
$\lim_{\b\to\infty}\frac{1}{\b}\log\frac{\mu_\b(\cS\backslash\cW)}{\mu_\b(\cS)}=
-\infty$.
\el

\bpr
We will prove that $\lim_{\b\to\infty}\frac{1}{\b}\log\mu_\b(\cS\backslash\cW)
/\mu_\b(\boxminus)=-\infty$. Since $\boxminus\in\cS$, this will prove the claim.

Let $f(\b)$ be the function satisfying (\ref{fprop}). We begin by noting that
\be{note}
\mu_\b(\cS\backslash\cW) \leq \mu_\b(\cI) \quad \mbox{ with } \quad
\cI = \big\{\s\in\cS\colon| \textnormal{supp}[C_B(\s)]|
> |\L_\b|/f(\b)\big\},
\ee
because the bootstrap percolation map increases the number of $(+1)$-spins. Let
$\cD(k)$ denote the set of configurations whose support consists on $k$ non-interacting
subcritical rectangles. Put $C_1=(\ell_c+2)(\ell_c+1)$. Since the union of a subcritical
rectangle and its exterior boundary has at most $C_1$ sites, it follows that in $\cI$ there
are at least $|\L_\b|/C_1f(\b)$ non-interacting rectangles. Thus, we have
\be{oje}
\mu_\b(\cI) \leq \sum_{k=\frac{|\L_\b|}
{C_1f(\b)}}^{K_{\textnormal{max}}} F(k)
\quad \mbox{ with } \quad
F(k) = \frac{1}{Z_\b} \sum_{ {\s\in\cX_\b\colon} \atop {C(\s)\in\cD(k)} }
e^{-\b\,H_\b(\s)},
\ee
where $K_{\max} \leq |\L_\b|$.

Next, note that
\be{Fkestpre}
F(k) \leq (2^{C_1})^k\,\frac{1}{Z_\b}\sum_{\s\in\cD(k)}\,e^{-\b H_\b(\s)}.
\ee
Since the bootstrap percolation map is downhill, the energy of a subcritical rectangle
is bounded below by $C_2=2J-h$ (recall Fig.~\ref{fig:gap}), and the number of ways to
place $k$ rectangles in $\L_\b$ is at most $\binom{|\L_\b|}{k}$, it follows that
\be{f1}
F(k) \leq 2^{C_1k}\,\binom{|\L_\b|}{k}\,\mu_\b(\boxminus)\,e^{-C_2\b k}
\leq 2^{C_1k}\,(C_1 e f(\b))^k\,\mu_\b(\boxminus)\,e^{-C_2\b k}
\leq \mu_\b(\boxminus)\,\exp[-\tfrac12 C_2\,\b k],
\ee
where the second inequality uses that $k! \geq k^k e^{-k}$, $k\in\N$, and the third
inequality uses that $f(\b)=e^{o(\b)}$. We thus have
\be{term1}
\sum_{k=\frac{|\L_\b|}{C_1f(\b)}}^{K_{\textnormal{max}}}\,F(k)\,\,
\leq 2\mu_\b(\boxminus)\,f(\b)\,\frac{|\L_\b|}{f(\b)}\,
\exp\left[-\tfrac12 \frac{C_2}{C_1}\,\b\,\frac{|\L_\b|}{f(\b)}\right],
\ee
which is the desired estimate because $|\L_\b|/f(\b)$ tends to infinity.
\epr


\section{Appendix: equivalence of ensembles and typicality of holes}
\label{appB}

\setcounter{equation}{0}

For $m\in\N$, let
\be{sm}
\cS^{(m)} = \big\{\s\in\cX_\b^{(m)}\colon\,|\textnormal{supp}[\s] \cap
B_{L_\b,L_\b}(x)|\leq \ell_c(\ell_c-1)+1\,\,\forall\,x\in\L_\b\big\}
\ee
and
\be{yset}
\begin{aligned}
\check{C}^{(m)}(0) &= \big\{\s\1_{\in B_{L_\b,L_\b}(0)}\colon\,\s\in\cS^{(m)}\big\},\\
\check{Z}_\b^{(m)}(0) &=\sum_{\s\in \check{C}^{(m)}(0)}
e^{-\beta\,H(\sigma)}.
\end{aligned}
\ee
The latter is the partition sum restricted to $B_{L_\b,L_\b}(0)$ when it carries $m$
particles. In Appendix~\ref{appB.1} we prove a lemma about ratios of partition sums
that was used in(\ref{muCnrellim}), (\ref{gamma2}), (\ref{I1u2}) and (\ref{Kawlbrels2alt}).
In Appendix~\ref{appB.2} we prove that $\lim_{\b\to\infty}\mu_\b(\check{\cS}(0))/\mu_\b
(\cS)=1$, which is needed in the proof of this lemma.

\subsection{Equivalence of ensembles}
\label{appB.1}

Recall (\ref{rhodef}), (\ref{cchec}) and (\ref{Cldef}).

\bl{cgcequiv}
$\check{Z}_\b^{(n_\b-s)}(0)/Z_\b^{(n_\b)} = (\rho_\b)^s\,\mu_\b(\cS)\,[1+o(1)]$ as
$\b\to\infty$ for all $s\in\N$.
\el

\bpr
The proof proceeds via upper and lower bounds.

\medskip\noindent
\underline{Upper bound:}
Let
\be{checkSdef}
\check{\cS}(0) = \big\{\s\in\cS^{}\colon\,\textnormal{supp}[\s] \cap B_{L_\b,L_\b}(0)
= \emptyset\big\}.
\ee
Write
\be{mucheckSrew}
\mu_\b(\check{\cS}(0)) = \frac{1}{Z_\b^{(n_\b)}} \sum_{\check{\s}\in\check{\cC}(0)}
\sum_{ {\zeta \subset [B_{L_\b,L_\b}(0)]^c\backslash\textnormal{supp}[\check{\s}]}
\atop {|\zeta|=s} }
\binom{n_\b}{s}^{-1}\,\1_{\{\check{\s}\vee\zeta\in\check{\cS}(0)\}}\,
e^{-\b H_\b(\check{\s}\vee\zeta)}.
\ee
This relation simply says that $n_\b$ particles can be placed outside $B_{L_\b,L_\b}(0)$
by first placing $n_\b-s$ particles and then placing another $s$ particles. Because the
interaction is attractive, we have
\be{Hattr}
H_\b(\check{\s}\vee\zeta) \leq H_\b(\check{\s}) + H_\b(\zeta)
\mbox{ and } H_\b(\zeta) \leq 0, \qquad \forall\,\check{\s},\zeta.
\ee
Consequently,
\be{mucheckSrewest}
\mu_\b(\check{\cS}(0)) \geq \binom{n_\b}{s}^{-1}\,\frac{1}{Z_\b^{(n_\b)}}
\sum_{\check{\s}\in\check{\cC}(0)} e^{-\b H_\b(\check{\s})}
\sum_{ {\zeta \subset [B_{L_\b,L_\b}(0)]^c\backslash\textnormal{supp}[\check{\s}]}
\atop {|\zeta|=s} }
\1_{\{\check{\s}\vee\zeta\in\check{\cS}(0)\}}.
\ee

We next estimate the second sum, uniformly in $\check{\s}$. When we add the $s$ particles,
we must make sure not to violate the requirement that all boxes $B_{L_\b,L_\b}(x)$,
$x\in\L_\b$, carry not more than $K$ particles (note that $\check{\cS}(0)\subset\cS$
and recall Definition~\ref{Rdef*}(a)). Partition $\L_\b\backslash B_{L_\b,L_\b}(0)$
into boxes of size $L_\b$. The total number of boxes containing $K$ particles is
at most $n_\b/K$. Hence, the total number of sites where we cannot place a particle
is at most $(n_\b/K)(3L_\b)^2$. Therefore
\be{sumestcom}
\sum_{ {\zeta \subset [B_{L_\b,L_\b}(0)]^c\backslash\{\check{\s}\}} \atop {|\zeta|=s} }
\1_{\{\check{\s}\vee\zeta\in\check{\cS}(0)\}}
\geq \binom{|\L_\b\backslash B_{L_\b,L_\b}(0)|-n_\b-(n_\b/K)(3L_\b)^2}{s},
\qquad \forall\,\check{\s}.
\ee
But $n_\b L_\b^2=o(n_\b/\rho_\b)=o(|\L_\b|)$ and $L_\b^2=o(1/\rho_\b)=o(|\L_\b|)$ by
(\ref{rhodef}) and (\ref{Lbetadef}--\ref{dbprop}), and so the right-hand side of
(\ref{sumestcom}) equals $[1+o(1)]\,|\L_\b|^s/s!$ as $\b\to\infty$. Since the binomial
in (\ref{mucheckSrewest}) equals $[1+o(1)]\,(n_\b)^s/s!$ with $n_\b=\lceil\rho_\b|\L_\b|
\rceil$, we end up with (recall (\ref{Cldef}))
\be{ubfinrev}
\mu_\b(\check{\cS}(0)) \geq \frac{\check{Z}_\b^{(n_\b-s)}(0)}{Z_\b^{(n_\b)}}
\,(\rho_\b)^{-s}\,[1+o(1)],
\ee
or
\be{ubfin}
\frac{\check{Z}_\b^{(n_\b-s)}(0)}{Z_\b^{(n_\b)}}
\leq (\rho_\b)^s\,\mu_\b(\check{\cS}(0))\,[1+o(1)].
\ee
Since $\check{\cS}(0)\subset\cS$, this gives the desired upper bound.

\medskip\noindent
\underline{Lower bound:}
Return to (\ref{mucheckSrew}). Among the $s$ particles that are added to $[B_{L_\b,L_\b}(0)]^c$,
let $s_1$ count the number that interact with the $n_\b-s$ particles already present and
$s_2$ the number that interact among themselves, where $s_1+s_2\leq s$. We can then estimate
\be{mucheckSrew3}
\begin{aligned}
&\mu_\b(\check{\cS}(0))\\
&\leq \frac{1}{Z_\b^{(n_\b)}}
\sum_{\check{\s}\in\check{\cC}(0)}\binom{n_\b}{s}^{-1}
e^{-\b H_\b(\check{\s})}\sum_{ {s_1,s_2} \atop {0 \leq s_1+s_2 \leq s} }
\bigg(\frac{s!}{s_1!\,s_2!}\bigg)^{-1}\,\\
&\qquad\qquad\times
\sum_{ {\zeta \subset [B_{L_\b,L_\b}(0)]^c\backslash\textnormal{supp}[\check{\s}]}
\atop {|\zeta|=s} } e^{-\b H(\zeta)}\,
\1_{\{|\zeta\cap\partial\check{\sigma}|=s_1\}}\,
\1_{\{s_2\,\,\textnormal{interacting particles in}\,\,\zeta\}}
\,\1_{\{\check{\s}\vee\zeta\in\check{\cS}(0)\}}\,\\
&\leq [1+o(1)]\,\frac{\check{Z}^{(n_\b-s)}_\b(0)}{Z_\b^{(n_\b)}}\,(\rho_\b)^{-s}\\
&\qquad\qquad
+ \frac{1}{Z_\b^{(n_\b)}} \sum_{\check{\s}\in\check{\cC}(0)}\binom{n_\b}{s}^{-1}
\,\,e^{-\b H_\b(\check{\s})}\sum_{ {s_1,s_2} \atop {1\leq s_1+s_2\leq s} }\\
&\qquad\qquad\times
\sum_{ {\zeta \subset [B_{L_\b,L_\b}(0)]^c\backslash\textnormal{supp}[\check{\s}]}
\atop {|\zeta|=s} } e^{-\b H(\zeta)}\,
\1_{\{|\zeta\cap\partial\check{\sigma}|=s_1\}}\,
\1_{\{s_2\,\,\textnormal{interacting particles in}\,\,\zeta\}}\,
\1_{\{\check{\s}\vee\zeta\in\check{\cS}(0)\}},
\end{aligned}
\ee
where in the second inequality we estimate the term with $s_1=s_2=0$ by using the
result for the upper bound. We will show that the other terms are exponentially
small.

For fixed $\check{\sigma}$, we may estimate the last sum in (\ref{mucheckSrew3}) as
\be{sumup}
\begin{aligned}
&\sum_{ {\zeta \subset [B_{L_\b,L_\b}(0)]^c\backslash\textnormal{supp}[\check{\s}]}
\atop {|\zeta|=s} } \,e^{-\b H(\zeta)}\,
\1_{\{|\zeta\cap\partial\check{\sigma}|=s_1\}}\,
\1_{\{s_2\,\,\textnormal{interacting particles in}\,\,\zeta\}}\,
\1_{\{\check{\s}\vee\zeta\in\check{\cS}(0)\}}\\
&\qquad\qquad\leq |\L_\b|^{s-s_1-s_2}\,(4n_\b)^{s_1}\,
\sum_{\sigma_\in \cS^{(s_2)}}
e^{-\b H(\sigma)}\,
\1_{\{s_2\,\,\textnormal{interacting particles in}\,\,\sigma\}}.
\end{aligned}
\ee
Indeed, $|\L_\b|^{s-s_1-s_2}$ bounds the number of ways to place $s-s_1-s_2$
non-interacting particles, and $(4n_\b)^{s_1}$ the number of ways to place $s_1$
particles that are interacting with the $n_\b-s$ particles already present.
Next, we write
\be{sumup2}
\begin{aligned}
&\sum_{\sigma_\in \cS^{(s_2)}}
e^{-\b H(\sigma)}\,
\1_{\{s_2\,\,\textnormal{interacting particles in}\,\,\sigma\}}\\
&\qquad\qquad
= \sum_{m=1}^{s_2}\sum_{j=1}^{m}\sum_{ {2\leq k_1,\ldots,k_j\leq K} \atop {\sum_{i=1}^j k_i=m} }
\sum_{ {C=\cup_{i=1}^j C_i} \atop {|C_i|=k_i\,\forall\,i}} e^{-\b\,\sum_{i=1}^j H(C_i) },
\end{aligned}
\ee
which is a cluster expansion of the partition function (with non-interacting clusters
$C_i$, all of which have size $\leq K=\ell_c(\ell_c+1)+1$). By a standard isoperimetric
inequality we have $H(\cC_i) \geq H_{k_i}$, with the latter denoting the energy of a droplet
of $k_i=|C_i|$ particles that is closest to a square or quasi-square. Hence
\be{sumup3}
\begin{aligned}
|\L_\b|^{-s_2}
&\sum_{\sigma_\in \cS^{(s_2)}}
e^{-\b H(\sigma)}\,
\1_{\{s_2\,\,\textnormal{interacting particles in}\,\,\sigma\}}\\
&\leq|\L_\b|^{-s_2}
\sum_{m=1}^{s_2}\sum_{j=1}^{m}\sum_{ {2 \leq k_1,\ldots,k_j\leq K} \atop {\sum_{i=1}^j k_i=s_2} }
e^{-\b\,\sum_{i=1}^j H_{k_i}}\,
\Bigg(\sum_{{C=\cup_{i=1}^j C_i} \atop { |C_i|=k_i\,\forall\,i} } \,1\Bigg)\\
&\leq\,C\,|\L_\b|^{-s_2}
\sum_{m=1}^{s_2}\sum_{j=1}^{m}\sum_{ {2\leq k_1,\ldots,k_j\leq K} \atop {\sum_{i=1}^j k_i=s_2} }
e^{-\b\,\sum_{i=1}^j H_{k_i}}\,|\L_\b|^j\\
&\leq C\,\sum_{m=1}^{s_2}\,\sum_{j=1}^{m}
\,\sum_{ {2 \leq k_1,\ldots,k_j\leq K} \atop {\sum_{i=1}^j k_i=s_2} }
e^{-\b\,\sum_{i=1}^j [H_{k_i}+(k_i-1)\b^{-1}\log|\L_\b|]}\\
&\leq \sum_{m=1}^{s_2}\,\sum_{j=1}^{m}
\,\sum_{ {2 \leq k_1,\ldots,k_j\leq K} \atop {\sum_{i=1}^j k_i=s_2} }
e^{-\b\,\sum_{i=1}^j [H_{k_i}+(k_i-1)\Delta]},
\end{aligned}
\ee
where in the last inequality we insert the bound $\b^{-1}\log|\L_\b|\geq \Delta$, which
is a immediate from (\ref{rhodef}) and (\ref{LambdacondsKa}).

Now, $H_{k_i}+k_i\Delta$ is the \emph{grand-canonical} energy of a square or quasi-square
with $k_i$ particles. It was shown in the proof of Proposition 2.4.2 in Bovier, den Hollander
and Nardi~\cite{BdHN06} that this energy is $\geq U\sqrt{k_i}$ for $1 \leq k_i \leq 4K$,
i.e., for a droplet twice the size of the proto-critical droplet. Since $2U>\Delta$, we
therefore have that $H_{k_i}+(k_i-1)\Delta>0$ when $k_i\geq 4$. Since $\Delta>U$, $H_2=-U$
and $H_3=-2U$, we have that also the terms with $k_i=2,3$ are $>0$. Consequently, there
exist $\epsilon>0$ and a constant $C$ that is independent of $\beta$ such that
\be{sumup5}
|\L_\b|^{-s_2}\sum_{\sigma_\in \cS^{(s_2)}}
e^{-\b H(\sigma)}\,
\1_{\{s_2\,\,\textnormal{interacting particles in}\,\,\sigma\}} e^{-\b H(\sigma)}
\leq C\,\,e^{-\b\,\epsilon}.
\ee
Combining (\ref{mucheckSrew3}--\ref{sumup}) and (\ref{sumup5}), we see that the correction
term in (\ref{mucheckSrew3}) is
\be{corr}
\begin{aligned}
&\mu_\b(\check{\cS}(0))-[1+o(1)]\,
\frac{\check{Z}^{(n_\b-s)}_\b(0)}{Z_\b^{(n_\b)}}\,(\rho_\b)^{-s}\\
&\qquad \leq C\,[1+o(1)]\,\frac{\check{Z}^{(n_\b-s)}_\b(0)}{Z_\b^{(n_\b)}}(\rho_\b)^{-s}
\sum_{{s_1,s_2}\atop{1\leq s_1+s_2\leq s}}(e^{U\b}\rho_\b)^{s_1}\,e^{-\b\epsilon}.
\end{aligned}
\ee
Since $\Delta>U$, we have $e^{U\b}\rho_\b\leq 1$ and so the sum is $o(1)$. Hence
\be{ZcheckZlb}
\frac{\check{Z}^{(n_\b-s)}_\b(0)}{Z_\b^{(n_\b)}}
\geq (\rho_\b)^s\,\mu_\b(\check{\cS}(0))\,[1+o(1)].
\ee
The claim now follows by using Lemma~\ref{checkSrel} below.
\epr

\subsection{Typicality of holes}
\label{appB.2}

\bl{checkSrel}
$\lim_{\b\to\infty}\mu_\b(\check{\cS}(0))/\mu_\b(\cS)=1$.
\el

\bpr
Since $\check{\cS}(0)\subset\cS$, we have $\mu_\b(\check{\cS}(0))\leq \mu_\b(\cS)$. It
therefore remains to prove the lower bound. Write
\be{lemb21}
\begin{aligned}
\mu_\b(\cS)
&= \mu_\b(\check{\cS}(0))\\
&\qquad +\sum_{m=1}^K \sum_{\eta\in\cX_\b^{(m)}}
\sum_{ {\zeta\in\cX_\b^{(n_\b-m)} } \atop {\eta\vee\zeta\in\cS} }
\frac{e^{-\b\,H(\eta\vee\zeta)}}{Z_\b^{(n_\b)}}\,
\1_{\{\textnormal{supp}[\eta]\subset B_{L_\b,L_\b}(0)\}}
\1_{\{\textnormal{supp}[\zeta]\subset [B_{L_\b,L_\b}(0)]^c\}}\\
&\leq \mu_\b(\check{\cS}(0))+ \g_1(\b) + \g_2(\b),
\end{aligned}
\ee
where
\be{g1errdef}
\g_1(\b) =\sum_{m=1}^K\sum_{\eta\in\cX_\b^{(m)}}
\sum_{ {\zeta\in\cX_\b^{(n-m)}} \atop {\eta\vee\zeta\in\cS} }
\frac{e^{-\b\,[H(\eta)+H(\zeta)]}}{Z_\b^{(n_\b)}}\,
\1_{\{\textnormal{supp}[\eta]\subset B_{L_\b,L_\b}(0)\}}
\1_{\{\textnormal{supp}[\zeta]\subset [B_{L_\b,L_\b}(0)]^c\}}
\ee
and $\g_2(\b)$ is a term that arises from particles interacting accross the
boundary of $B_{L_\b,L_\b}(0)$. We will show that both $\g_1(\b)$ and $\g_2(\b)$
are negligible.

Estimate, with the help of (\ref{ubfin}) (and recalling (\ref{sm}--\ref{yset})),
\be{lemb22}
\begin{aligned}
\g_1(\b) &\leq \sum_{m=1}^K \frac{\check{Z}_{\b}^{(n_\b-m)}}{Z_\b^{(n_\b)}}
\sum_{\eta\in\cS^{(m)}} e^{-\b H(\eta)}\,
\1_{\{\textnormal{supp}[\eta]\subset B_{L_\b,L_\b}(0)\}}\\
&= [1+o(1)]\,\mu_\b(\check{\cS}(0))
\sum_{m=1}^K (\rho_\b)^m \sum_{\eta\in\cS^{(m)}} e^{-\b\,H(\eta)}\,
\1_{\{\textnormal{supp}[\eta]\subset B_{L_\b,L_\b}(0)\}}\\
&= [1+o(1)]\,\mu_\b(\check{\cS}(0))
\sum_{m=1}^K (\rho_\b)^m \sum_{j=1}^m
\sum_{ {2\leq k_1,\ldots,k_j\leq K} \atop {\sum_{i=1}^j k_i=m} }
\sum_{ {C=\cup_{i=1}^j C_i} \atop {|C_i|=k_i\,\forall\,i} } e^{-\b\,\sum_{i=1}^j H(C_i)},
\end{aligned}
\ee
where the last equality is a cluster expansion as in (\ref{sumup2}). Using once more
the isoperimetric inequality, we get (recall (\ref{Lbetadef}))
\be{lemb24}
\begin{aligned}
\g_1(\b) &\leq
[1+o(1)]\,\mu_\b(\check{\cS}(0))\,
\sum_{m=1}^K (\rho_\b)^m \sum_{j=1}^m
\sum_{ {2\leq k_1,\ldots,k_j\leq K} \atop {\sum_{i=1}^j k_i=m} }
e^{-\b\sum_{i=1}^j H(k_i)}
\Bigg(\sum_{ {C=\cup_{i=1}^j C_i} \atop {|C_i|=k_i\,\forall\,i} }\,1\Bigg)\\
&\leq C\,\mu_\b(\check{\cS}(0))
\sum\in_{m=1}^K (\rho_\b)^m \sum_{j=1}^m (L_\b^2)^j\,
\sum_{ {2\leq k_1,\ldots,k_j\leq K} \atop {\sum_{i=1}^j k_i=m} }
e^{-\b\sum_{i=1}^j H_{k_i}}\\
&= C\,\mu_\b(\check{\cS}(0)) \sum_{m=1}^K\,\sum_{j=1}^m
\sum_{ {2\leq k_1,\ldots,k_j\leq K} \atop {\sum_{i=1}^j k_i=m} }
e^{-\b\sum_{i=1}^j [H_{k_i}+k_i\Delta-(\Delta-\delta_\b)]}\\
&\leq C'\,\mu_\b(\check{\cS}(0))\,e^{-\b\epsilon}
\end{aligned}
\ee
for some $\epsilon>0$ and constants $C,C'<\infty$ that are independent of $\b$.

Estimate, with the help of (\ref{ubfin}),
\be{lemb26}
\begin{aligned}
\g_2(\b) &\leq \sum_{m=1}^K \sum_{\eta\in\cS^{(m)}}
e^{-\b H(\eta)} \sum_{k=1}^m e^{\b kU}\,
\1_{\{\textnormal{supp}[\eta]\subset B_{L_\b,L_\b}(0)\}}\,\,
\frac{\check{Z}_{\b}^{(n_\b-(m+k))}}{Z_\b^{(n_\b)}}\\
&\leq \sum_{m=1}^K \sum_{\eta\in\cS^{(m)}}
e^{-\b H(\eta)} \sum_{k=1}^m e^{\b kU}\,
\1_{\{\textnormal{supp}[\eta]\subset B_{L_\b,L_\b}(0)\}}\,\,
(\rho_\b)^{m+k}\,\mu_\b(\check{\cS}(0))\,[1+o(1)]\\
&\leq [1+o(1)]\,\mu_\b(\check{\cS}(0))\sum_{m=1}^K (\rho_\b)^m
\sum_{\eta\in\cS^{(m)}} e^{-\b H(\eta)}
\sum_{k=1}^m e^{-\b k(\Delta-U)}
\1_{\{\textnormal{supp}[\eta]\subset B_{L_\b,L_\b}(0)\}},
\end{aligned}
\ee
and we can proceed as (\ref{lemb22}--\ref{lemb24}) to show that this term is negligible.
\epr

\subsection{Atypicality of critical droplets}
\label{appB.3}

The following lemma was used in Section~\ref{S4.1.1}.

\bl{GibbscC+}
$\lim_{\b\to\infty} \mu_\b(\tilde{\cC}\backslash\cC^+)/\mu_\b(\cS)=0$.
\el

\bpr
Similarly as in (\ref{lemb21}), we first write
\be{Ctildests}
\begin{aligned}
&\mu_\b(\tilde{\cC}\backslash\cC^+) \leq \mu_\b(\tilde{\cC})\\
&= |\L_\b|\,\g(\b) + |\L_\b| \sum_{\eta\in\cX_\b^{(K)}}
\sum_{ {\zeta\in\cX_\b^{(n_\b-K)}} \atop {\eta\vee\zeta\in\tilde{\cC}} }
\frac{e^{-\b\,[H(\eta)+H(\zeta)]}}{Z_\b^{(n_\b)}}\,
\1_{\{\textnormal{supp}[\eta]\subset B_{L_\b,L_\b}(0)\}}
\1_{\{\textnormal{supp}[\zeta]\subset [B_{L_\b,L_\b}(0)]^c\}}.
\end{aligned}
\ee
with $\g(\b)$ a negligible error term that arises from particles interacting accross
the boundary of $B_{L_\b,L_\b}(0)$. We then proceed as in (\ref{g1errdef}--\ref{lemb24}),
obtaining ($\Gamma=\Gamma^*+K\Delta$)
\be{finestproc}
\begin{aligned}
\mbox{r.h.s.(\ref{Ctildests})}
&\leq N\,|\L_\b|\,e^{-\b\Gamma^*}\,(\rho_\b)^K\,
\mu_\b(\check{\cS}(0))\,[1+o(1)]\\
&= N\,|\L_\b|\,e^{-\b\Gamma}\,\mu_\b(\cS)\,[1+o(1)],
\qquad \b\to\infty,
\end{aligned}
\ee
which is $o(\mu_\b(\cS))$ by (\ref{LambdacondsKa}).
\epr


\section{Appendix: Typicality of starting configurations}
\label{appC}

\setcounter{equation}{0}

In Sections~\ref{appC.1}--\ref{appC.2} we prove the claims made in the remarks
below (\ref{L*def}), respectively, (\ref{L*defKa}).

\subsection{Glauber}
\label{appC.1}

\bpr
Split
\be{ssplit}
\cS=\cS_L\cup(\cS\setminus\cS_L)=\cS_L\cup\cU_{>L},
\ee
where $\cU_{>L}\subset\cS$ are those configurations $\s$ for which $C_B(\s)$ has at
least one rectangle that is larger than $Q_L(0)$. We have
\be{cboots}
C_B(\s) = \bigcup_{x\in X(\s)} R_{\ell_1(x),\ell_2(x)}(x),
\ee
where $X(\s)$ is the set of lower-left corners of the rectangles in $C_B(\s)$, which in
turn can be split as
\be{Xsplit}
X(\s) = X^{>L}(\s) \cup X^{\leq L}(\s),
\ee
where $X^{>L}(\s)$ labels the rectangles that are larger than $Q_L(0)$ and $X^{\leq L}
(\s)$ labels the rest.

Let $\s|_{A}$ denote the restriction of $\s$ to the set $A\subset\Z^2$. Then, for any
$x\in X(\s)$, we have
\be{Hsumsplit}
H(\s) = H\big(\s|_{R_{\ell_1(x),\ell_2(x)}(x)}\big)
+ H\big(\s|_{R^c_{\ell_2(x),\ell_2(x)}(x)}\big),
\ee
because the rectangles in $C_B(\s)$ are non-interacting. Since for $\s\in\cU_{>L}$ there
is at least one rectangle with lower-left corner in $X^{>L}(\s)$, we have
\be{ulsplit}
\begin{aligned}
\mu_\b(\cU_{>L})
&\leq \sum_{x\in\L_\b} \sum_{\s\in\cS}
\1_{\{x\in X^{>L}(\s)\}}\,\mu_\b(\s)\\
&= \sum_{x\in\L_\b} \sum_{\s\in\cS}
\1_{\{x\in X^{>L}(\s)\}}\,
\frac{1}{Z_\b}\,\exp\Big\{-\b\Big[H\big(\s|_{R_{\ell_1(x),\ell_2(x)}(x)}\big)
+ H\big(\s|_{R^c_{\ell_1(x),\ell_2(x)}(x)}\big)\Big]\Big\}\\
&\leq e^{-\b\Gamma_{L+1}} \sum_{x\in\L_\b} \sum_{\s\in\cS}
\1_{\{x\in X^{>L}(\s)\}}\,\frac{1}{Z_\b}
e^{-\b H\big(\s|_{R^c_{\ell_1(x),\ell_2(x)}(x)}\big)},
\end{aligned}
\ee
where $\Gamma_{L+1}$ is the energy of $Q_{L+1}(0)$. In the last step we use the fact that
the bootstrap map is downhill and that the energy of $Q_L(0)$ is increasing with $L$.
Since the energy of a subcritical rectangle is non-negative, we get
\be{ulsplit2}
\mu_\b(\cU_{>L}) \leq N_{L+1}\,e^{-\b\Gamma_{L+1}}\,|\L_\b|\,\mu_\b(\cS)
\ee
with $N_{L+1}$ counting the number of configurations with support in $Q_{L+1}(0)$.

On the other hand, by considering only those configurations in $\cU_{>L}$ that have a
$Q_{L+1}(0)$ droplet, we get
\be{lowul}
\mu_\b(\cU_{>L}) \geq N_{L+1}\,e^{-\b\Gamma_{L+1}}\,|\L_\b|\,\mu_\b^{[Q_{L+1}(0)]^c}(\cS),
\ee
where the last factor is the Gibbs weight of the configurations in $\cS$ with support
outside $[Q_{L+1}(0)]^c$. It easy to show that $\mu_\b^{[Q_{L+1}(0)]^c}(\cS)=\mu_\b(\cS)
[1+o(1)]$ as $\b\to\infty$ and so
\be{ulsplit4}
\mu_\b(\cU_{>L}) \geq N_{L+1}\,e^{-\b\Gamma_{L+1}}\,|\L_\b|\,\mu_\b(\cS)\,[1+o(1)],
\qquad \b\to\infty.
\ee
Combining (\ref{ulsplit2}--\ref{lowul}), we conclude that $\lim_{\b\to\infty}
\mu_\b(\cU_{>L})/\mu_\b(\cS)=0$ if and only if
\be{volcond}
\lim_{\b\to\infty} |\L_\b|\,e^{-\Gamma_{L+1}}=0.
\ee
\epr

\subsection{Kawasaki}
\label{appC.2}

\bpr
Split
\be{ssplitalt}
\cS=\cS_L\cup(\cS\setminus\cS_L)=\cS_L\cup\cU_{>L},
\ee
where $\cU_{>L}\subset\cS$ are those configurations $\s$ for which there exists
an $x$ such that $|{\rm supp}[\s]\cap B_{L_\b,L_\b}(x)|>L$. Then
\be{ul2split}
\mu_\b(\cU_{>L})
\leq \sum_{x\in\L_\b} \sum_{\s\in\cS} \sum_{m=L+1}^{K}
\mu_\b(\s)\,\1_{\{|{\rm supp}[\s]\cap B_{L_\b,L_\b}(x)|=m\}}
=|\L_\b|\,[\varphi(\b)+\g(\b)],
\ee
where
\be{phideful}
\varphi(\b) = \sum_{m=L+1}^K \sum_{\eta\in\cX_\b^{(m)}}
 \sum_{ {\zeta\in\cX_\b^{(n_\b-m)}} \atop {\eta\vee\zeta\in\cS} }
\frac{e^{-\b [H(\eta)+H(\zeta)]}}{Z_\b^{(n_\b)}}\,
\1_{\{\textnormal{supp}[\eta]\subset B_{L_\b,L_\b}(0)\}}\,
\1_{\{\textnormal{supp}[\zeta]\subset [B_{L_\b,L_\b}(0)]^c\}}
\ee
and $\g(\b)$ is an error term arising from particles interacting accross
the boundary of $B_{L_\b,L_\b}(0)$. By the same argument as in (\ref{lemb26}),
this term is negligible. Moreover,
\be{ul3split}
\begin{aligned}
\varphi(\b)
&\leq \sum_{m=L+1}^K \frac{\check{Z}_\b^{(n_\b-m)}}{Z_\b^{(n_\b)}}
\bigg(\sum_{\eta\in\cS^{(m)}}
e^{-\b\,H(\eta)}\,
\1_{\{\textnormal{supp}[\eta]\subset B_{L_\b,L_\b}(0)\}}\bigg)\\
&\leq  [1+o(1)]\,\mu_\b(\cS) \sum_{m=L+1}^K (\rho_\b)^m\,
\bigg(\sum_{\eta\in\cS^{(m)}} e^{-\b H(\eta)}\,
\1_{\{\textnormal{supp}[\eta]\subset B_{L_\b,L_\b}(0)\}}\bigg),
\end{aligned}
\ee
where in the last inequality we use Lemmas~\ref{cgcequiv}--\ref{checkSrel}.
Now proceed as in (\ref{lemb22}--\ref{lemb24}), via the cluster expansion,
to get
\be{ul5split}
\begin{aligned}
\varphi(\b)
&\leq 1+o(1)]\,C\,\mu(\cS) \sum_{m=L+1}^{K} \sum_{j=1}^m
\sum_{{2 \leq k_1,\ldots,k_j \leq K} \atop {\sum_{i=1}^j k_i=m} }
e^{-\b[H_{k_i}+k_i\Delta-(\Delta-\delta_\b)]}\\
&\leq [1+o(1)]\,C\,\mu(\cS)\, e^{-\b[\Gamma_{L+1}-(\Delta-\delta_\b)]},
\end{aligned}
\ee
where $H_k$ is the energy of a droplet with $k$ particles that is closest to
a square or quasi-square, $\Gamma_{L+1}= H_{L+1}+(L+1)\Delta$, and the second
inequality uses the isoperimetric inequality together with the fact that
$H_k+k\Delta$ is increasing in $k$ for subcritical droplets.

On the other hand, by considering only those configurations in $\cU_{>L}$ that
have a droplet with $L+1$ paticles, we get
\be{ul6splitalt}
\varphi(\b) \geq [1+o(1)]\,C\,\mu(\cS)\, e^{-\b[\Gamma_{L+1}-(\Delta-\delta_\b)]}.
\ee
Combining (\ref{ul2split}) and (\ref{ul5split}--\ref{ul6splitalt}), we conclude
that $\lim_{\b\to\infty}\mu_\b(\cU_{>L})/\mu_\b(\cS)=0$ if and only if
\be{volcondalt}
\lim_{\b\to\infty} |\L_\b|\,e^{-\b\,(\Gamma_{L+1}-(\Delta-\delta_\b))}=0.
\ee
\epr


\section{Appendix: The critical droplet is the threshold}
\label{appD}

\setcounter{equation}{0}

In this appendix we show that our estimates on capacities imply that the \emph{average}
probability under the Gibbs measure $\mu_\b$ of destroying a supercritical droplet and
returning to a configuration in $\cS_L$ is exponentially small in $\b$. We will give the
proof for Kawasaki dynamics, the proof for Glauber dynamics being simpler.

Pick $M\geq\ell_c$. Recall from (\ref{eqm.2}) that $e_{\cD_M,\cS_L}(\s)=c_\b(\s)\P_\s
\left(\t_{\cS_L}<\t_{\cD_M}\right)$ for $\s\in\cD_M$. Hence summing over $\s\in
\partial\cD_M$, the internal boundary of $\cD_M$, we get using (\ref{eqm.3}) that
\be{over.2}
\frac{\sum_{\s\in\partial \cD_M}\mu_\b(\s)
c_\b(\s)\P_\s\left(\t_{\cS_L}<\t_{\cD_M}\right)}
{\sum_{\s\in\partial\cD_M} \mu_\b(\s)c_\b(\s)}
=\frac{\CAPA(\cS_L,\cD_M)}
{\sum_{\s\in\partial\cD_M} \mu_\b(\s)c_\b(\s)}.
\ee
Clearly, the left-hand side of (\ref{over.2}) is the escape probability to $\cS_L$ from
$\partial\cD_M$ \emph{averaged} with respect to the canonical Gibbs measure $\mu_\b$
conditioned on $\partial\cD_M$ weighted by the outgoing rate $c_\b$. To show that this
quantity is small, it suffices to show that the denominator is large compared to
the numerator.

By Lemma~\ref{CAPSRcest*},
\be{over.1}
\CAPA(\cS_L,\cD_M) \leq \CAPA(\cS_L,(\cS^c\setminus\tilde{\cC})\cup\cC^+)
= N\,|\L_\b|\,\frac{4\pi}{\D\b}\,e^{-\b\Gamma}\,\mu_\b(\cS)[1+o(1)].
\ee
On the other hand, note that $\partial\cD_M$ contains all configurations $\s$ for which
there is an $M \times M$ droplet somewhere in $\L_\b$, all $L_\b$-boxes not containing
this droplet carry at most $K$ particles, and there is a free particle somewhere in
$\L_\b$. The last condition ensures that $c_\b(\s)\geq 1$. Therefore we can use
Lemma~\ref{cgcequiv} to estimate
\be{over.3}
\sum_{\s\in\cD_M} \mu_\b(\s)c_\b(\s) \geq |\L_\b|\,e^{-\b H_{M^2}}\,
\frac{\check{Z}_\b^{(n_\b-M^2)}}{Z_\b^{(n_\b)}}
= |\L_\b|\,e^{-\b H_{M^2}}\,
(\rho_\b)^{M^2}\,\mu_\b(\cS)\,[1+o(1)],
\ee
where $H_{M^2}$ is the energy of an $M \times M$ droplet. Combining (\ref{over.1}--\ref{over.3})
we find that the left-hand side of (\ref{over.2}) is bounded from above by
\be{over.4}
\left(N\frac{4\pi}{\D\b}\right)\,\frac{\exp\left[-\b \G\right]}
{\exp\left[-\b(H_{M^2}+\Delta M^2)\right]}\,[1+o(1)],
\ee
which is exponentially small in $\b$ because $\G>H_{M^2}+\Delta M^2$ for all $M\geq\ell_c$.


\end{document}